\documentclass[preprint,showpacs,showkeys,preprintnumbers,amsmath,amssymb]{revtex4}
\usepackage{latexsym}
\usepackage{amssymb}
\usepackage{amsfonts}
\usepackage{bm}

\usepackage{graphicx,epsfig}

\newtheorem{defi}{Definition}

\newtheorem{lem}{Lemma}
\newtheorem{theo}{Theorem}

\newtheorem{rem}{Remark}

\newcommand{\Xla}{X_{\lambda}}
\newcommand{\Xo}{X^{\ast}}
\newcommand{\noi}{\noindent}

\newcommand{\DP}{\displaystyle}

\newcommand{\vth}{L_{n}}

\newcommand{\de}{\delta}

\newcommand{\la}{\lambda}

\newcommand{\te}{\theta}

\newcommand{\om}{\omega}
\newcommand{\Om}{\Omega}

\newcommand{\ka}{\mathbb{K}}
\newcommand{\km}{k_{c}}

\newcommand{\R}{\mathbb{R}}
\newcommand{\E}{\mathbb{E}}

\newcommand{\vrc}{\mathbb{V}ar}
\newcommand{\cvc}{\mathbb{C}ov}

\def\bigone{\mbox {\large $1\hspace{-3.5pt} \mbox{l}$}}
\newcommand{\bfs}{{\bf s}}

\newcommand{\bfr}{{\bf r}}
\newcommand{\bfU}{{\bf U}}
\newcommand{\bfu}{{\bf u}}
\newcommand{\bfom}{{\bm \omega}}
\newcommand{\bftheta}{{\bm \theta}}

\newcommand{\rmx}{{\rm x}}

\newcommand{\bfz}{{ \hat{\bm \omega}}}
\newcommand{\dd}{d/2-1}
\def\czero{\textrm{$c_d^{(0)}$}}
\def\cone{\textrm{$c_d^{(1)}$}}
\def\ctwo{\textrm{$c_d^{(2)}$}}
\def\cthree{\textrm{$c_d^{(3)}$}}

\def\sigx{\sigma_{\rmx}^2}
\newcommand{\e}{{\eta_1}}
\newcommand{\fbar}{\overline{f_{X}}(h)}
\newcommand{\fbarg}{\overline{f_{X}}(h_1)}
\newcommand{\fbarc}{\overline{f_{X}}(h_2)}
\newcommand{\fbarcc}{\overline{f_{X}}(\sqrt{2}h_2)}
\newcommand{\fbarccc}{\overline{f_{X}}(2h_2)}
\newcommand{\fbarcr}{\overline{f_{X}}(rh_2)}
\newcommand{\phiso}{\overline{\varphi_1}(h_1)}
\newcommand{\phieo}{\phi_1(a_1)}
\newcommand{\phist}{\overline{\varphi_2}(h_2)}

\begin{document}

\title{On the Inference of Spatial Continuity using Spartan Random Field Models}

\author{Samuel Elogne}
\email{elogne@mred.tuc.gr }
\affiliation{Department of Mineral Resources Engineering\\
Technical University of Crete\\Chania 73100, Greece}%
\author{Dionissios T. Hristopulos}
\email{dionisi@mred.tuc.gr}
 \homepage{http://www.mred.tuc.gr/home/hristopoulos/dionisi.html}
\affiliation{Department of Mineral Resources Engineering\\
Technical University of Crete\\Chania 73100, Greece}%
\thanks{}

\begin{abstract}
This paper addresses the inference of spatial dependence in the
context of a recently proposed framework. More specifically, the
paper focuses on the estimation of model parameters for a class of
generalized Gibbs random fields \citep{dth03}, i.e., Spartan Spatial
Random Fields (SSRFs). The problem of parameter inference is based
on the minimization of a distance metric. The latter involves a
specifically designed distance between sample constraints (variance,
generalized ``gradient'' and ``curvature'') and their ensemble
counterparts. The general principles used in the construction of the
metric are discussed and intuitively motivated. In order to enable
calculation of the metric from sample data, estimators for
generalized ``gradient'' and ``curvature'' constraints are
constructed.  These estimators, which are not restricted to SSRFs,
are formulated using compactly supported kernel functions. An
intuitive method for kernel bandwidth selection is proposed. It is
proved that the estimators are asymptotically unbiased and
consistent for differentiable random fields, under specified
regularity conditions. For continuous but non-differentiable random
fields, it is shown that the estimators are asymptotically
consistent. The bias is calculated explicitly for different kernel
functions. The performance of the sample constraint estimators and
the SSRF inference process are investigated by means of numerical
simulations.
\end{abstract}

\keywords{non-parametric, inverse problem, optimization, simulation}


\maketitle

\section{Introduction}
Spatial Random Fields (SRF's) have a wide range of applications in
subsurface hydrology \citep{gel93,kit97,rub03}, petroleum
engineering \citep{hohn}, environmental data analysis
\citep{christ,smith00,kan04}, mining exploration and reserves
estimation \citep{goov,arm98}, and environmental health \citep{ch98}
among other fields. From the applications viewpoint, the main goals
are first to characterize the spatial continuity of such processes,
and then to exploit the continuity for spatial estimation
(prediction) and simulation. A methodological problem of continuing
interest is the inference of the random field parameters that
characterize the spatial continuity from the experimental data. The
latter are typically distributed on irregular sampling grids.

This paper seeks to address the issue of random field inference
within the context of a specific model, the
Fluctuation-Gradient-Curvature (FGC) Spartan Spatial Random Field
(SSRF), which was introduced in \cite{dth03}. SSRFs result from a
convolution of a kernel function with an underlying SRF that may
include non-resolved fine-scale detail at length scales below
$\lambda$. The kernel function acts as a low-pass filter that
suppresses the spectral component of fluctuations above a cutoff
wavevector $\km.$ The removed part corresponds to sub-resolution
scales. We will denote {\it sample averages} with a horizontal bar
over the averaged quantity, and {\it ensemble averages} with the
mathematical expectation symbol, i.e., $m_{\rm X}= \E[\Xla({\bf
s})]$. Without loss of generality in the following it is assumed
that $m_{\rm X}=0.$

Let the data be given by the measurements ${\bf \Xo} \equiv
\{\Xo_{1}, \ldots, \Xo_{n}\}$, of the scalar quantity $X$ at the set
of sampling locations ${\bf S}_{\rm m} \equiv \{\bfs_{1},\ldots
\bfs_{n}\}$ in the domain $\Omega_n \in \R^{d}$. The area enclosed
by the convex hull of $\Omega_n$ is denoted by $|\Omega_n|$. We
assume that the data can be modeled as a \textit{sample}
(realization) of ${X}_{\lambda}({\bf s})$, which is a Gaussian,
weakly stationary, isotropic FGC-SSRF. The isotropic assumption is
not a major restriction, since
 under certain conditions the anisotropic
parameters can be established and isotropic conditions can be
restored by rotation and rescaling transformations
\citep{dth02,dth04,dth05a}. The isotropic  FGC-SSRF  involves the
\textit{parameter set} ${\bm \theta}=(\eta_0, \eta_1, \xi, \km)$:
the \textit{scale coefficient} $\eta_0$ determines the variance, the
\textit{shape coefficient} $\eta_1$ determines the shape of the
covariance function, the \textit{characteristic length} $\xi$
determines the range of spatial dependence, and $\km$ is a
\textit{wavevector cutoff} related to the resolution scale $\lambda$
\cite{dth03}.

Regarding parameter inference, the main idea introduced in
\cite{dth03} and further elaborated here, is that the SSRF model
parameters can be estimated by matching sample constraints and their
ensemble counterparts.  The model parameters are determined by
treating the sample constraints as estimators of the respective
stochastic constraints. This perspective relies on the validity of
the \textit{ergodic hypothesis} \citep{yagl87,lantu}.

The constraint matching idea is similar to the standard approach, in
which the experimental variogram is matched with various model
functions to determine an optimal model of spatial continuity.
However, there are significant differences between the two
approaches. (1) In the SSRF approach the number of estimated
constraints is small (four in the case of the FGC-SSRF). This is
 due to the efficient parametrization of spatial dependence in the
 FGC-SSRF, which is
 based on interactions instead of the covariance matrix.
In contrast, variogram modeling attempts to match the entire
functional dependence of the variogram function.   (2) The FGC-SSRF
includes a family of covariance functions that account for various
types of spatial continuity \citep{dth03,dthse06}. Hence, in
practice fitting the sample constraints with one SSRF model may be
sufficient. In contrast, the experimental variogram is fitted with a
number of model functions to determine the ``optimal'' spatial
model. (3) The SSRF sample constraints focus on the short-range
behavior of spatial continuity. This is motivated by two
observations: in geostatistical applications, the long-range
behavior can only be estimated with significant uncertainty; in
addition, it is known that the long-range behavior does not have a
significant impact on optimal linear prediction in regions where the
field is densely sampled \citep{stein99}. (4) The computational
complexity of SSRF constraint calculations, at least on regular
grids, is $O(mn)$, where $m$ is $o(n)$ and depends on the kernel
bandwidth, while for variogram calculations the respective
complexity is at best $O(n\log n)$ if tree-based structures are
used, or $O(n^{2})$ using standard methods \cite{dth03}.

The main results obtained in this paper include the following: (1)
Generalized gradient and curvature estimators are formulated in
terms of kernel averages, and a consistent method for selecting the
kernel bandwidths is proposed.  The generalized gradient and
curvature estimators have a wider scope than the FGC-SSRF model:
They are defined for both differentiable and continuous (but
non-differentiable) spatial models. In the differentiable case, the
estimators are defined in terms of finite-difference approximations
of the respective derivatives. In the continuous case, the finite
differences are not divided by the corresponding length spacing
(step) in order to obtain asymptotically well defined quantities.
(2) Convergence properties for the generalized gradient and
curvature estimators are proved. (3) The constraint-based parameter
inference procedure introduced in \cite{dth03} is improved by adding
a constraint that eliminates the nonlinear dependence of the model
variance on the FGC-SSRF parameters. (4) Numerical simulations
establish the performance of the constraints estimators and the
parameter inference procedure.

The remaining of this paper is structured as follows: An
introduction to the FGC-SSRF in continuum space is given in
Section~(\ref{sec:fgc-ssrf}).  In Section~(\ref{sec:constraint}) the
definition of the sample constraints on regular grids is reviewed.
This is followed by the definition of generalized stochastic
constraints for the FGC-SSRF in Section~(\ref{sec:cons-irreg}).
Generalized sample constraints for the gradient and curvature are
defined in Section~(\ref{sec:sample}). Theorems establishing the
convergence of the constraint estimators are stated and proved in
Section~(\ref{sec:asympt}). Subsequently, the parameter inference
process developed in \cite{dth03} is reviewed and refined in
Section~(\ref{sec:inference}). Finally, numerical simulations are
used to validate the estimators and the parameter inference process
in Section~(\ref{sec:simul}).

\section{Review of the FGC-SSRF Model}
\label{sec:fgc-ssrf}
 In general, a \textit{Gibbs} random field has
the following joint probability density function (p.d.f.)
\begin{equation}
\label{gibbspdf} f_{\rm x} [X_\lambda({\bf s});{\bm \theta}] = \frac
        {\exp \left\{ { - H[X_\lambda ({\bf s});{\bm \theta}]} \right\} }
        {Z({\bm \theta})},
\end{equation}
where  $ H[X_\lambda({\bf s});{\bm \theta}] $ is the \textit{energy
functional},  ${\bm \theta}$ is a set of \textit{model parameters},
and the constant $ Z({\bm \theta}) $, called the \textit{partition
function} is the p.d.f. normalization factor obtained by integrating
$ \exp \left\{ -H[X_\lambda({\bf s});{\bm \theta}] \right\} $ over
all the realizations of the SRF. The FGC p.d.f. in $ R^d $ is
determined from the equation:

\begin{equation}
\label{fgc} H_{\rm fgc}[X_\lambda ({\bf s});{\bm \theta} ] =
\frac{1}{{2\eta _0 \xi ^d }} \int d{\bf s} \, \textsl{h}_{\rm fgc}
\left[ {X_\lambda({\bf s});{\bm \theta'} } \right],
\end{equation}
where  ${\bm \theta'}=(\eta _1 ,\xi,\km)$, and $ \textsl{h}_{\rm
fgc}$ is the normalized (to $\eta_0=1$) local energy density at
point $\bf{s}$. The functional $\textsl{h}_{\rm fgc} \left[
{X_\lambda({\bf s});{\bm \theta'} } \right] $ is given \textit{in
the continuum} by the following expression

\begin{equation}
\label{fgccont} h_{{\rm fgc}} \left[ {X_\lambda ({\bf s});{\bm
\theta'} } \right] = \left[ {X_\lambda ({\bf s})} \right]^2 + \eta
_1 \,\xi ^2 \left[ {\nabla X_\lambda ({\bf s})} \right]^2  + \xi ^4
\left[ {\nabla ^2 X_\lambda ({\bf s})} \right]^2 ,
\end{equation}
The functional (\ref{fgccont}) is permissible if the resulting
covariance function is positive definite, i.e., if it satisfies
\textit{Bochner's theorem} \citep{bochner}. Permissibility
constrains the value of $\e$ (see \citep{dth03},\citep{dthse06}).

The explicit, albeit non-linear, dependence of the p.d.f. on three
physically meaningful parameters, $\eta_0, \eta_1, \xi,$ instead of
three linear coefficients multiplying the terms $\left[ {X_\lambda
({\bf s})} \right]^2 ,$ $\left[ {\nabla X_\lambda ({\bf s})}
\right]^2$ and $\left[ {\nabla ^2 X_\lambda ({\bf s})} \right]^2,$
simplifies the parameter inference problem and allows intuitive
initial guesses for the parameters.

The FGC model has a particularly simple expression in Fourier space.
If the Fourier transform of the covariance function is defined by
means of
\begin{equation}
    \label{eq:covft0}
    \tilde{G}_{{\rm x};\lambda}({\bf k};{\bm \theta})= \int d{\bf r} e^{-\imath \bf{k\cdot r}}
    G_{{\rm x};\lambda}({\bf r};{\bm \theta}),
\end{equation}
then the energy functional in Fourier space is given by

\begin{equation}
\label{enfunc} H_{\rm fgc}[X_\lambda ({\bf s}) ;{\bm \theta} ] =
\frac{1}{2 (2 \pi)^{d}} \int d{\bf k} \tilde{X}_{\lambda}({\bf k})
\, \tilde{G}^{-1}_{{\rm x};\lambda}({\bf k};{\bm \theta}) \,
\tilde{X}_{\lambda}(-{\bf k}).
\end{equation}
The interaction is diagonal in Fourier space, i.e., the
\textit{precision matrix} $ \tilde{G}^{-1}_{{\rm x};\lambda}({\bf
k};{\bm \theta})$ couples only components with equal wavevectors.
For a real-valued SSRF $\Xla({\bf s})$ it follows that
$\tilde{X}_{\lambda}(-{\bf k})=\tilde{X}^{\dagger}_{\lambda}({\bf
k})$. Since Bochner's theorem guarantees the non-negativity of the
covariance spectral density, it follows from~(\ref{enfunc}) that the
energy is a non-negative functional.

The covariance spectral density follows from the expression:

\begin{equation}
\label{covspd} \tilde G_{{\rm x;}\lambda}({\bf k};{\bm \theta}) =
\frac{{\left| {\,\tilde Q_\lambda({\bf k})\,} \right|^{{\kern 1pt}
2} \,\eta _0 \,\xi ^d }}{{1 + \eta _1 \,(k{\kern 1pt} \xi )^2  +
(k{\kern 1pt} \xi )^4 }}
\end{equation}
where  $ \tilde Q_\lambda ({\bf{k}})\ $ is the Fourier transform of
the coarse-graining kernel, which determines how the fluctuations
are cut off at the resolution scale $\lambda$  \citep{dth03}. For
isotropic SSRF's,  $\tilde Q_\lambda ({\bf k})$ has no directional
dependence. In \citep{dth03,dth05b}, a kernel having a boxcar
spectral density with a sharp wavevector cutoff  at $\km $ was used.
This kernel leads to a band-limited covariance spectral density  $
\tilde G_{{\rm x;}\lambda }({\bf k};{\bm \theta}) $.

\section{Generalized Energy Functional}
\label{sec:constraint}

The energy density defined by Eq.~(\ref{fgccont}) is valid in the
continuum, and for differentiable SRFs. Generalized versions of the
functional that are valid  on regular lattices can be defined. For
example:
\begin{defi}\label{def1}
 We define the local energy terms
$S_0,$ $S_1(a_1),$ and $S_2(a_2)$, $\forall \bfs \in
\mathbb{L}_{d}$, as follows:
$$
S_0=\left[ \Xla({\bf s})\right]^2, \quad S_1(a_1)= \sum_{i=1}^d \Big
[ \Xla({\bf s}+a_1 {\bf \vec{e}}_i)-
   \Xla({\bf s}) \Big ]^2/a_1^2 $$
$$ S_2(a_2)= \sum_{i,j=1}^d\Delta_2^{(i)}[\Xla({\bf s})] \;
\Delta_2^{(j)}[\Xla({\bf s})],
$$
where $\Delta_2^{(i)}$ is the centered second-order difference
operator in the direction ${\bf \vec{e}}_i,$ i.e,
\begin{equation*}\label{deltai}\Delta_2^{(i)}
[\Xla({\bf s})]=\Big [\Xla({\bf s}+a_2\;{\bf\vec{e}}_i)
        +\Xla({\bf s}-a_2\;{\bf\vec{e}}_i)
        -2\Xla({\bf s})\Big ]/a_2^2.
\end{equation*}
 $S_0$ represents the square of the fluctuations, $S_1$ the
square of the generalized gradient, and $S_2$ the square of the
generalized curvature.
\end{defi}

The \textit{generalized gradient} and \textit{curvature} terms above
are expressed in terms of finite differences instead of derivatives.
These terms replace the gradient and curvature in (\ref{fgccont}).
On a hypercubic lattice $\mathbb{L}_{d} \subset\mathbb{Z}^{d}$ in
$d$ dimensions with step $a$, one obtains $a_1=a_2=a$. The sample
counterparts of $S_{i}(\bfs), \, i=0,1,2$, obtained by replacing
$\Xla({\bf s})$ with $\Xo(\bfs)$, are thus well-defined even for
non-differentiable SRFs. Parameter inference is based on matching
the sample constraints, $\overline{\mathcal{S}_{i}(\bfs)}$, with the
stochastic constraints, $\E[S_{i}(\bfs)]$, as shown in
\citep{dth03}.

\begin{defi}\label{def2}
The ensemble moments  $\E[S_0],$   $\E[S_1(a_1)]$ and $\E[S_2(a_2)]$
provide the SSRF model constraints. These can be expressed in terms
of the variance $G_{\la}(0)$ and the semivariogram function
$F_{\la}$  as follows:

\begin{equation}\E[S_0]= G_{\la}(0),\qquad
 \label{sto1}
 \end{equation}
 \begin{equation}
 \E[S_1(a_1)]:=\frac{\phi_1(a_1)}{a_1^{2}} = \frac{\cone}{a_1^{2}}\;  F_{\la}(a_1)
 \label{sto2}
 \end{equation}
\begin{equation}
 \label{sto3}
\E[S_2(a_2)] :=\frac{\phi_2(a_2)}{a_2^{4}}= \frac{1}{a_2^4} \,
\left[\ctwo \, F_{\la}(a_2) - \cthree \, F_{\la}(\sqrt{2} a_2)-
\cone\, F_{\la}(2a_2) \right],
 \end{equation}
where $\cone=2d,$ $\ctwo=8d^2$, and $\cthree=4d(d-1).$
\end{defi}

\begin{rem} {\rm The SSRF constraints are expressed in terms of the
semivariogram $F_{\la}$, but this does not imply that the
experimental variogram is required for determining the spatial
model}.
\end{rem}

The dependence of the stochastic constraints on the SSRF parameters
$\bftheta$ is not shown explicitly to keep the notation concise. The
stochastic constraints are well defined for the FGC-SSRF, which are
differentiable if $\km$ is finite \cite{dthse06}. In the case of
continuous but non-differentiable models, the ratios
$\phi_1(a_1)/a_{1}^{2}$ and $\phi_2(a_2)/a_{2}^{4}$ are not well
defined when $a_1,a_2\rightarrow 0$.  Then, the constraints are
defined in terms of the quantities $\phi_1(a_1)$ and $\phi_2(a_2)$
respectively.

\section{Constraint Definitions on Irregular Grids}
\label{sec:cons-irreg}

In most geostatistical applications the available sample is
distributed on an irregular sampling grid. In order to infer the
model parameters, suitable stochastic and sample-based constraints
need to be defined. On regular grids the lattice symmetry leads to
obvious choices for the {\it distance increments} (steps) $a_1$ and
$a_2$ and the finite differences.  On irregular grids, we formulate
the sample gradient and curvature constraints using kernel averages.
We also define steps $a_1$ and $a_2$, suitable for general sampling
point distributions. In addition, the kernel bandwidths are selected
so as to yield good asymptotic properties for the generalized
gradient and curvature estimators.

\subsection{Stochastic FGC-SSRF Constraints}
\label{ssec:stoch-con} The stochastic constraints are related to the
SRF model and thus do not depend on the sampling point distribution.
Hence, the constraints defined in (\ref{sto1})-(\ref{sto3}) can be
used for irregular grids as well. The dependence of the constraints
on the SSRF parameters is made explicit using the \textit{spectral
representation of the covariance function} \citep{yagl87}. If we
define $v=(\km\xi)^{2}$, $\Pi(v)=1+\eta_1 v+ v^2,$ and
$w=v^{1/2}\xi^{-1}$ then for any $a> 0$ the following relations
hold:
\begin{equation}
\label{eq:cov}G_{\la}(a ) = \frac{\eta_0\, \xi^{\dd}
}{2\,(2\pi)^{d/2}\,a^{\dd}} \,
  \int_{0}^{\infty }dv\, v^{(d-2)/4} \, J_{\dd}(a\, w)
  \,\frac{ \big |\tilde{Q}_{\la}(w) \, \big
    |^2}{\Pi(v) }.
\end{equation}
The \textit{variance stochastic constraint}, obtained for $a=0$, is
given by
\begin{equation}
\label{stoch-0} \E[S_0]=\frac{\eta_0 }{2^d\,\pi^{d/2} \Gamma(d/2)}
\,
  \int_{0}^{\infty }dv\, v^{\dd}
  \, \frac{ \big |\tilde{Q}_{\la}(w) \, \big
    |^2}{\Pi(v)}.
\end{equation}
In general, the dependence of $\E[S_0]$ on $\e$ and $\km\xi$ is
nonlinear \citep{dthse06}, i.e., $\sigx= \eta_0 \, g_{d}(\e,\km
\xi)$. The function $g_{d}(\e,\km \xi)$ tends to an asymptotic
finite bound as $\km \xi \rightarrow \infty$. The bound is attained
with very good accuracy if $\km \xi = q_d$, where $q_d$ is an $O(1)$
constant that depends on $d$. To eliminate the dependence of the
SSRF variance on $\e$ and $\km\xi$, we impose the relation

\begin{equation}
\label{eq:eta0}\eta_0=2^d\,\pi^{d/2} \, \Gamma(d/2)\, \sigx,
\end{equation}
 which is equivalent to the following normalization constraint:
\begin{equation}\label{s0prime}
\E[S_0^{'}]=\int_{0}^{\infty }dv\, v^{\dd}
  \,  \frac{ \big |\tilde{Q}_{\la}(w) \, \big|^2}{\Pi(v)}=1.
\end{equation}

 The  stochastic constraint for the \textit{generalized gradient} is given by
\begin{equation}
\label{semi-der2} \E[S_1(a_1)] =\frac{\cone\sigx } {a_{1}^{2}}
\int_{0}^{\infty }dv\, \left [ v^{\dd}-\gamma_d \, v^{(d-2)/4} \,
\frac{J_{\dd}(a_1 w)} {a_1^{\dd}}\right]\frac{\big
|\tilde{Q}_{\la}(w) \, \big
    |^2}{\Pi(v)},
\end{equation}
where $\gamma_d= (2\xi)^{\dd}\Gamma(d/2)$.

Based on (\ref{sto3}) and (\ref{eq:cov}), the stochastic constraint
for the \textit{generalized curvature} is given by
\begin{eqnarray}
\label{semi-der4}  \E[S_2(a_2)] & = & \frac{\sigx } {a_{2}^{4}}
\int_{0}^{\infty }dv\, \frac{\big |\tilde{Q}_{\la}(w) \, \big
|^2}{\Pi(v)} \, \left\{ \left[\cthree+3\cone \right]
v^{\dd}-\gamma_d\, \frac{v^{(d-2)/4}}{a_2^{\dd}} \right. \nonumber\\
    & &
    \left. \left[ \ctwo J_{\dd}(a_2w)   - \cthree \, \frac{J_{\dd}(a_2 \sqrt{2}w)}{2^{d/4-1/2}}
    - \cone\, \frac{J_{\dd}(2a_2\, w)}{2^{\dd}}  \right] \right\}.
\end{eqnarray}

The selection of $\eta_0$ based on Eq.~(\ref{eq:eta0}) increases the
number of SSRF constraints to four; the other three are given by the
equations~~(\ref{stoch-0}), (\ref{semi-der2}) and~(\ref{semi-der4}).
Thus, the number of parameters matches the number of constraints.

The steps $a_1$ and $a_2$ depend on the sampling point distribution.
Their selection is discussed in Section~(\ref{sec:asympt}) below. In
general, $a_2$ is different from $a_1$. To incorporate the spatial
modelling of data from non-differentiable distributions, one should
focus on the quantities $\phi_1(a_1)=a_1^{2} \, \E[S_1(a_1)]$ and
$\phi_2(a_2)=a_2^{4} \, \E[S_2(a_2)]$ instead of $\E[S_1(a_1)]$ and
$\E[S_2(a_2)].$

\section{Sample Constraints}
\label{sec:sample} We formulate sample constraints that provide
`well-behaved' estimators of the model constraints defined above. We
emphasize that the following sample estimators of the generalized
gradient and curvature constraints are not restricted to the
FGC-SSRF model.

The variance constraint is local, i.e., it does not involve
differences between neighboring points. Hence, if the distribution
of the sampling points is uniform, it is sufficient to use the
classical variance estimator. To estimate a non-zero mean one can
use the sample average, $\hat{m}_{\rm x}=n^{-1}\sum_{i}^{n} \Xo({\bf
s}_{i})$, in view of which the sample variance
$\overline{\mathcal{S}_0}$ is given by:
\begin{equation}
\label{barS0} \overline{\mathcal{S}_0}=\frac{1}{n}\sum_{i=1}^n
\left[ \Xo({\bfs}_{i})-\hat{m}_{\rm x}\right]^2.
\end{equation}

\noindent Declustered estimates of the mean and the variance can be
used if the sampling point distribution is non-uniform, in order to
obtain unbiased estimates of the variance. However, non-ergodic
fluctuations, which often dominate the estimation of the variance
from a single sample, are not significantly reduced by cell
declustering. Another possibility is using the kernel-based variance
estimator \cite{hall1}, which has improved convergence properties.
However, we are not aware of a systematic method for selecting the
kernel bandwidth for the variance.

\subsection{Kernel Averages of Sample Functions}
\label{ssec:kern-av} To define sample-based generalized gradient and
curvature constraints we use \textit{isotropic kernel functions}
$K({\bf r})$ with suitably selected bandwidth parameters, $h_{1}$
(for the gradient estimator) and $h_{2}$ (for the curvature
estimator). The selection of the bandwidths is guided by consistency
principles that link them to the respective steps, $a_1$ and $a_2$.

The kernel $K$ is a bounded, positive, and compactly supported
$[0,R]$ function. Hence, the  moments
\begin{equation} \label{eq:mom1} m_{K,j}=
\int_0^{R} ds\;s^{j-1} \, K(s), \end{equation} and \begin{equation}
\label{eq:mom2} m^{(2)}_{K,j}= \int_0^{R} ds\;s^{j-1} \, K^{2}(s).
\end{equation}
are finite
 for all $j \in \mathbb{Z}^{+}$.  In addition, we define the
 {\it kernel moment ratio:}
 \begin{equation}
 \label{eq:mom-rat}
 B_p:= \frac{m_{K,d+p}}{m_{K,d}}.
 \end{equation}
In practice non-compactly supported kernels (e.g., the Gaussian
kernel,) that decrease to $0$ faster than polynomially work just as
well as compactly supported kernels.

The following notation is introduced to facilitate calculations with
kernel averages. For $\bfs_i, \bfs_j \in {\bf S}_{\rm m}$,
$\bfs_{i,j} := \bfs_i - \bfs_j$ will denote the distance vector, and
$s_{i,j}:= \|\bfs_{i,j}\|$ its Euclidean norm.

The abbreviations $W_{i,j} :=K \big ( {\bfs_{i,j}}/h_1\big ),$ and
$Q_{i,j\, (r)} :=K \big ( {\bfs_{i,j}}/rh_2\big )$ where
$r=1,2,\sqrt{2}$ will be used for the kernel weights. The weights
$W_{i,j}$ and $Q_{i,j\, (r)}$ are random variables, due to the
variability in the sampling positions.

The symbol $\sum'_{i,j}$ will denote a summation over both position
indices $i$ and $j$ excluding the diagonal terms $i=j$. Similarly,
the triple summation $\sum'_{i,j,k}$ and the quadruple summation
$\sum'_{i,j,k,l}$ will exclude all the terms in which at least two
indices take the same value.

Given kernel bandwidths $h_1$ and $h_2$, and a two-point sample
function $A_{i,j}=A(\Xo_{i},\Xo_{j})$  or a function of sampling
positions $A_{i,j}=A(\bfs_i, \bfs_j)$, where $i,j=1\ldots,n$, the
{\em off-diagonal kernel-weighted averages} will be denoted by:

\begin{equation}
\label{Ksum} \ka_{h_1} \left\{A_{i,j} \right\} :={\sum}'_{i, j}
W_{i,j}\, A_{i,j},
\end{equation}

\begin{equation}
\label{Ksum2} \ka_{r h_2} \left\{A_{i,j} \right\} :={\sum}'_{i, j}
Q_{i,j\,(r)}\, A_{i,j}.
\end{equation}
The \textit{normalized kernel average of  $A_{i,j}$}
 is defined by means of the equation:
 \begin{equation}
 \label{kernel-av}
\left\langle  A_{i,j} \right\rangle_{h_1} := \frac
 {\mathbb{K}_{h}\left\{ A_{i,j} \right\} }{\mathbb{K}_{h}\{1\}}.
 \end{equation}

More specifically, we will denote the \textit{sample increment SRF}
by means of
$$X^{*}_{i,j} := \Xo({\bf s}_i)-\Xo({\bf s}_j).$$
\noi The sample function that represents the kernel average of
$X^{*}_{i,j}$ with a kernel bandwidth $h$ will be denoted by
\begin{equation}
\label{fbar} \fbar:=\left\langle \, [\Xo_{i,j}]^{2} \,
\right\rangle_{h}.
\end{equation}
The random variable $\fbar$ incorporates variability due to both
$\Xo$ and the sampling positions.

\subsection{Definition of Sample Gradient and Curvature Estimators}
\label{ssec:sample-cons}

The notation introduced above is now used to define sample
estimators for the squares of the generalized gradient and
curvature.

\begin{defi}\label{def-gradcons}
The sample average of the square generalized gradient is defined as
follows:

\begin{equation}
\label{barS1irr} \overline{\mathcal{S}_1}(a_1):=
\frac{\cone}{2a_{1}^{2}} \,
 \left\langle  \left [X^{*}_{i,j}\right]^2
 \right\rangle_{h_1}=\frac{\cone}{2a_{1}^{2}} \,
 \fbarg:=\frac{\overline{\varphi_1}(h_1)}{a_{1}^{2}},
 \end{equation}
 where $\overline{\varphi_1}(h_1)=d\,\fbarg$.

 The bandwidth $h_1$ is related
 to $a_1$ by means of the following
 consistency principle:
\begin{equation}
\label{eq:alpha1} a_1^{2} =\left\langle s_{i,j}^2
\right\rangle_{h_1}.
\end{equation}
 \end{defi}

The step-bandwidth dependence introduced by the consistency
principle is physically motivated, because only $a_1$ represents an
actual length scale. By adopting (\ref{eq:alpha1}), the kernel
average in (\ref{barS1irr}) is forced to focus on points separated
by distances controlled by the step value. This makes sense for
calculations of generalized gradient and curvature terms. The sample
constraint defined in (\ref{barS1irr}) is analogous to the
respective stochastic constraint in (\ref{sto2}). In addition, the
consistency principle ensures that, for differentiable SRFs,
$\overline{\mathcal{S}_1}(a_1)$ is an unbiased estimator of
$\E[S_1(a_1)]$ (see below).

We introduced the three related quantities, $\fbarg$,
$\overline{\varphi_1}(h_1)$ and $\overline{\mathcal{S}_1}(a_1)$ for
the following reasons: If the observed SRF can be considered
differentiable, the generalized gradient constraint
$\overline{\mathcal{S}_1}(a_1)$ is well defined at the limit $a_1
\rightarrow 0.$  If the observed SRF is continuous but
non-differentiable the limit of $\overline{\mathcal{S}_1}(a_1)$ as
$a_1 \rightarrow 0$ does not exist. In this case, it makes more
sense to work with the kernel-averaged square increment,
$\overline{\varphi_1}(h_1)$.  To simplify the accounting it is often
advantageous to work with the square increment per direction,
$\fbarg$; in the isotropic case $\overline{\varphi_1}(h_1)$ and
$\fbarg$ are simply proportionally related. Similar comments apply
to the case of the generalized curvature constraint.
\\

\begin{defi}\label{def-curvcons}
The sample average of the  square generalized curvature is defined
as follows:
\begin{equation}
\label{barS2irr}
\overline{\mathcal{S}_{2}}(a_2):=\frac{\overline{\varphi_2}(h_2)}{a_{2}^{4}}=
\frac{1}{2a_{2}^{4}} \left[ \ctwo \mu_1(h_2)  \,
 \fbarc-\cthree \mu_2(h_2)
 \fbarcc-\cone \fbarccc \right],
 \end{equation}

\noi the constants $\mu_1(h_2), \mu_2(h_2)$ are given by the
following averages:

\begin{equation} \label{muhat1}
\mu_1(h_2)= \frac{\DP \cthree \mu_2(h_2) \left\langle s_{i,j}^2
\right\rangle_{\sqrt{2}h_2}+\cone\left\langle s_{i,j}^2
\right\rangle_{2h_2}} {\DP \ctwo \left\langle s_{i,j}^2
\right\rangle_{h_2}},
\end{equation}
\begin{equation} \label{muhat2}
\mu_2(h_2)= \frac{\DP [\ctwo+8\cone]\left\langle s_{i,j}^4
\right\rangle_{h_2}+\cone\left\langle s_{i,j}^4
\right\rangle_{h_2}\frac{\left\langle s_{i,j}^2
\right\rangle_{2h_2}}{\left\langle s_{i,j}^2
\right\rangle_{h_2}}-\cone\left\langle s_{i,j}^4
\right\rangle_{2h_2}} {\DP \cthree \left\langle s_{i,j}^4
\right\rangle_{\sqrt{2}h_2}-\cthree \left\langle s_{i,j}^4
\right\rangle_{h_2}\frac{\left\langle s_{i,j}^2
\right\rangle_{\sqrt{2}h_2}}{\left\langle s_{i,j}^2
\right\rangle_{h_2}}}.
\end{equation}

The bandwidth $h_2$ is linked to the step by means of the
consistency principle:
\begin{equation}\label{alpha2}
a_2^{4} = \left\langle s_{i,j}^4 \right\rangle_{h_2},
\end{equation}
\end{defi}

\noi The sample constraint $\overline{\mathcal{S}_{2}}(a_2)$ given
by (\ref{barS2irr}) includes three terms that correspond to the
terms in the respective stochastic constraint (\ref{sto3}). The
coefficients $\mu_1(h_2)$ and  $\mu_2(h_2)$ in (\ref{barS2irr})
incorporate the impact of the sampling network topology and the
kernel function used. As shown in Lemma~(\ref{lem:mu1mu2-as}),
Section~(\ref{ssec:gen-curv}), the coefficients $\mu_1(h_2)$ and
$\mu_2(h_2)$ are approximately equal to $1$. Their precise form is
selected to ensure that, in the case of differentiable SRFs, the
generalized curvature constraint is asymptotically unbiased.

\section{Asymptotic Properties of Constraint Estimators}
\label{sec:asympt}

\vspace{0.5cm}

The asymptotic limit corresponds to $n \rightarrow \infty$. At the
limit it is assumed that $1/|\Om_n| \rightarrow 0$. In order to
establish the asymptotic properties of the sample estimators for the
generalized gradient and curvature, we first present some formalism
and the conditions required for the validity of the proofs.
\\

\subsection{Formalism} \label{notation} The following notation will
be used in the proofs of asymptotic behavior.

\begin{enumerate}

\item \label{c:loc} The sampling
locations will be expressed as ${\bf s}_i = \vth {\bf u}_i$, where
$\vth \propto |\Omega_n|^{1/d}$ is the characteristic domain scale,
and $\bfu_i$ denote the realizations of the random vector $\bfU_i
\supset [0,1]^{d}$. \\

\item \label{c:dist} For any vectors ${\bf v}_{i}$, ${\bf v}_{j}$,
 the pair distance will be denoted by
  ${\bf v}_{i,j}:= {\bf v}_i - {\bf v}_j$,  and its Euclidean
norm by $v_{i,j}:= \| {\bf v}_i - {\bf v}_j \|$.
\\

\item \label{c:uij} In integrals of kernel averages,
the distance of the normalized sampling locations will be
denoted by $\bfom:=\bfu_i - \bfu_j$. \\

\item \label{c:int}
A vector $\bfom$ will be expressed in spherical polar coordinates as
$$\bfom=\om \,\bfz\quad \mbox{and}\quad (\bfz)_i=\cos\te_i \, {\prod}_{0\leq j<i}\sin\te_j $$
where $\te_0=\pi/2,$ $\te_d=0,$ $\te_{d-1}\in [0,2\pi)$, and $ \te_i
\in [0,\pi)$ for $i=1,\ldots,d-2.$ The Jacobian of the
transformation is given by $\om^{d-1}J_d(\bftheta)$ where
$\bftheta=\big (\te_1,\ldots,\te_{d-1} \big )$ and
$$J_d(\bftheta)=\left (\sin\,\te_1\right)^{d-2}\left (\sin\,\te_2\right)^{d-3}
\ldots \sin\,\te_{d-2}.$$ The  area of the $d-$dimensional unit
sphere will be denoted by:
$$  \quad A_d:=\int_{\mathbb{S}_d}d\bftheta
J_d(\bftheta),$$ where
$\int_{\mathbb{S}_d}:=\int_0^{\pi}\ldots\int_0^{\pi}
\int_{0}^{2\pi}. $
\\

\item \label{c:pert}
The following aspect ratios will be used as small
\textit{perturbation parameters}: $p_n:=h_1/\vth$ and $q_n :=
h_2/\vth$.
\\

\end{enumerate}

\subsection{Conditions} The following conditions will be assumed
to hold:

\begin{enumerate}
\item \label{H:loc}
The normalized location random vectors ${\bf U}_1, \ldots, {\bf
U}_n$ are assumed to be independent and identically distributed.
\\

\item \label{H:f1}
  The probability density function (pdf)
  $f_1({\bf u}_{i,j})$ of the sampling-location pair-distance vector
is continuously differentiable in a neighborhood of zero.\\

\item \label{H:f2}
The joint pdfs $f_2({\bf u}_{i,j},{\bf u}_{i,k})$ and $f_3({\bf
u}_{i,j},{\bf u}_{i,k},{\bf u}_{i,l})$ are also continuously
 differentiable in a neighborhood of zero. \\

\item  \label{H:gc}
The conditional pdf $g_{\bfu}(\bfu_{i,j} |\bfu_j=\bfu)$
 is uniformly bounded in $\bfu.$\\

\item  \label{H:match}
The joint moments of $\Xo$ are identical to those of $X_{\la}.$ For
example, $\forall \,\, \bfs_{i}, \bfs_{j}$, with $i \neq j$,
$F^{\ast}(\bfs_{i} - \bfs_{j})= F_{\la}(\bfs_{i} - \bfs_{j}). $ \\

\item  \label{H:semi}
The model semivariogram $F_{\la}$ is continuous
in a neighborhood of zero.\\

\item  \label{H:psil}
There exists a continuous and bounded function $\psi_{\la}$ such
that \begin{equation} \label{psila}\mathbb{C}ov  \left[
(\Xo_{i,j})^2, \, (\Xo_{p,q})^2 \right] =
 \psi_{\la}\left(s_{i,p},s_{j,q},s_{i,q},s_{j,p}\right ).
 \end{equation}

  For example, if  $\Xla$ is a Gaussian SRF with
semivariogram $F_{\la},$

\begin{equation*}\label{psi}
    \psi_{\la}(u_1,u_2,u_3,u_4 )=2\Big[ F_{\la}(u_1)+
 F_{\la}(u_2)-F_{\la}(u_3)-
 F_{\la}(u_4)\Big]^2.
\end{equation*}

\item \label{H:psil-prop}
If $\epsilon<<1$,
 there exist $c_1\geq 1$
 and three continuous functions $g_1,$ $g_2$ and $g_3$
 such that:
$$\psi_{\la}(0,0,\epsilon\,s,\epsilon\,s)=g_1(s) \left (\frac{\epsilon}{\xi}\right
)^{2c_1}+o\left (\frac{\epsilon}{\xi}\right )^{2c_1},$$
$$\psi_{\la}(\epsilon s_1,\epsilon s_2,0,\epsilon s_3)=g_2(s_1,s_2,s_3) \left
(\frac{\epsilon}{\xi}\right )^{2c_1}+o\left
(\frac{\epsilon}{\xi}\right )^{2c_1}
$$ and
\begin{eqnarray*}& &\psi_{\la}\left ( u_3,
\big\|\epsilon \bfom_2+u_3\bfz_3- \epsilon \bfom_1\big\|,
  \big\|u_3\bfz_3+  \epsilon \bfom_2\big\|,
   \big\|u_3\bfz_3-\epsilon \bfom_1\big\|\right )=\\
& & g_3(u_3,\om_1,\om_2,\bftheta_1,\bftheta_2,\bftheta_3) \, \left
(\frac{\epsilon}{\xi}\right )^{2c_1}+o\left
(\frac{\epsilon}{\xi}\right )^{2c_1}.
\end{eqnarray*}
 For example, if  $\Xla$ is a Gaussian SRF  with a spherical or
exponential covariance function, the above conditions hold with
$c_1=1.$  In the case of differentiable covariance models commonly
used (Gaussian, hole-type, rational quadratic, Cauchy) one has
$c_1=2.$ \\

\item \label{H:g3bounds}
The following integral of the function $g_{3}$ is bounded:

 \item \label{H:band}
The bandwidths $h_1$ and $  h_2 $ tend to $0$ as
$n$ tends to $\infty$.\\

\item \label{H:dens} At the asymptotic limit,
$  | \Om_n| /n $, tends to $0$ as $n$ tends to $\infty$. This
condition is satisfied simultaneously with $1/|\Om_n | \rightarrow
0,$  if $|\Om_n | \propto n^{\de}$
and $0 <\de <1.$\\

 \end{enumerate}

Conditions (\ref{H:loc})-(\ref{H:gc}) specify properties of the
sampling point distribution. Condition (\ref{H:match}) expresses the
correspondence between the SRF model and the sampled data.
Conditions (\ref{H:semi})-(\ref{H:g3bounds}) specify properties
which are satisfied by default for FGC-SSRFs. They are explicitly
stated here, because the convergence properties of the constraint
estimators are proved for more general cases, including non-Gaussian
SRF models. In particular, conditions
(\ref{H:psil})-(\ref{H:g3bounds}) are used in the analysis of the
sample constraints variance. Conditions
(\ref{H:band})-(\ref{H:dens}) imply an {\em asymptotic
densification} of the sampling network, since the area enclosed by
the convex hull increases slower than the number of points. For
regular grids, this condition is obtained if the spacing decreases
as the number of nodes increases. The densification conditions are
necessary for proving asymptotic convergence of the estimators.

\begin{lem}
\label{lem:ergod} If the above conditions hold, the following is
true:

\begin{equation}
\label{eq:sumW}
  {\rm Pr}\Big( \lim_{n\rightarrow \infty} \ka_{h_1}\{A_{i,j}\}
 = n^{2} \, \E\left[W_{i,j} \, A_{i,j}\right] \Big) = 1,
\end{equation}
where the indices $i,j$ refer to any pair of non-identical sampling
points. Similarly, if the summation is over a weighted $k-point$ $(k
\in \mathbb{Z})$ non-diagonal function, the result is $\propto
n^{k}$.
 \end{lem}
This ergodic result follows directly by applying
 the arguments in the proof of Theorem 3.1 in \citet{hall1},
 which will not be repeated here. Equation~(\ref{eq:sumW}) enables
 the calculation of sample kernel averages in the asymptotic limit.

 We will also use the following lemma:
\begin{lem}
\label{lem:useful} Let $X_n$ be a sequence of uniformly bounded
random variables such that $X_n=o(1)$ almost surely. Then
$\E[X_n]^k=o(1),$ for any $k.$
\end{lem}

\subsection{Generalized Gradient Estimation}
\label{ssec:gen-grad} In this section we prove a relation between
the ``gradient'' step and the kernel bandwidth $h_1$, and we propose
a physical estimate for the step. We also investigate the asymptotic
properties of the mean and variance of the generalized gradient
estimator.

\begin{lem}\label{alpha1-as}
 The following relation holds between the bandwidth $h_1$ and the gradient
 step $a_1$:
\begin{equation}\label{eq:bandw1}
 a_1=h_1 B_2^{1/2}+O(h_1\,p_n)
 \quad {\rm a.s.},
 \end{equation}
 where $B_2$, is the kernel moment ratio
defined in (\ref{eq:mom-rat}).

 \end{lem}

{\bf Proof:}

Based on the consistency principle (\ref{eq:alpha1}), the step $a_1$
is expressed as follows:

\begin{equation}
\label{eq:Kha1-W} a_1^{2}= \frac{ \sum'_{i, j} W_{i,j} \, s_{i,j}^2}
     {  \sum'_{i, j} W_{i,j} }.
\end{equation}

The above can be calculated explicitly in the asymptotic regime
using Lemma~(\ref{lem:ergod}).

 \vspace{6pt}

\textit{Leading-order calculation of  $\E[W_{i,j}]$.} \vspace{4pt}
$$ \E[W_{i,j}] =  \int d\bfom \;
  K\big (\vth\|\bfom \|/h_1\big) \, f_1(\bfom).$$

\noi The dominant asymptotic contribution from the integral is
evaluated by means of a Taylor expansion of the pdf $f_1$ in terms
of the small parameter $p_n$ using the condition (\ref{H:f1}), i.e.,

\begin{eqnarray}
\label{eq:ewijs}
  \E[W_{i,j}] &=&  \int_0^{\infty} d\om \;\om^{d-1} K\left( \om /p_n\right)
 \int_{\mathbb{S}_d}d\bftheta\;
  J_d(\bftheta)f_1\left(\om \bfz \right)
  \nonumber \\
 &=&p_n^d \int_0^{R} du \;u^{d-1} K\left( u\right)
 \int_{\mathbb{S}_d}d\bftheta\;
  J_d(\bftheta)f_1\left(p_n\,u \bfz \right)
  \nonumber \\
  & = & p_n^{d} \, A_d \, f_1(0) \, m_{K,d}
  +O\left(p_n^{d+1}\right).
\end{eqnarray}

\noindent This expansion gives the asymptotically dominant term of
$\E[W_{i,j}]$. Based on Lemma~(\ref{lem:ergod}) it follows that:

\begin{equation}
\label{eq:ewij}
  {\sum_{i,j}}' \, W_{i,j}=n^2 \, p_n^{d} \, A_d \, f_1(0) \, m_{K,d}
  +O\left(n^2p_n^{d+1}\right) \quad \mbox{a.s.},
\end{equation}
where $A_d$ is the area of the $d-$dimensional unit sphere defined
in paragraph~(\ref{c:int}) of the Notation subsection.

\textit{Leading-order calculation of
$\E\left[W_{i,j}\,s_{i,j}^{2}\right]$.} \vspace{4pt}

\begin{eqnarray}
\label{eq:esswijs}
 \E\left[W_{i,j}\, s_{i,j}^{2} \right] &=&
 \vth^2 \, \E\left[ s_{i,j}^{2} \,
 K\left (p_n  \, U_{i,j} \right )
   \right]  = \vth^2 \, \int d\bfom \;\|\bfom\|^2
  K\left(\|\bfom\| /p_n\right) \,  f_1(\bfom)  \nonumber \\
  &=& \vth^2  p_{n}^{d+2} A_d f_1(0) \, m_{K,d+2}
  +O\left(\vth^2 \, p_{n}^{d+3} \right).
\end{eqnarray}
Hence, from Lemma~(\ref{lem:ergod}) and equation (\ref{eq:esswijs})
it follows that
\begin{equation}
\label{eq:esswij}
  {\sum_{i,j}}' W_{i,j}\, s_{i,j}^{2}=n^2\vth^2  p_{n}^{d+2} A_d f_1(0) \, m_{K,d+2}
  +O\left(n^2\vth^2 \, p_{n}^{d+3} \right)\quad \mbox{a.s.}
\end{equation}
Based on (\ref{eq:Kha1-W}), (\ref{eq:ewij}), and (\ref{eq:esswij})
the asymptotic behavior of the step $a_1$ is given by:

\begin{equation*}\label{eq:asym-h1}
    a_1^2=h_1^2 \, B_2
    +O(h_1^2p_n),\quad \mbox{a.s.}.
\end{equation*}

 \vspace{4pt}

\subsubsection{Selection of Distance Step}
\label{sssec:step}
 Lemma (\ref{alpha1-as}) is valid for any step
$a_1$. We define by $\mathfrak{B}_0$ the set that includes for every
sampling point $\bfs_i$ the distance vectors from all its near
neighbors $\bfs_j$, and also $N_0=\big |\mathfrak{B}_0\big |.$ A
sensible estimate $\hat{a}_1$ is the geostatistical $d$-power
average of the Euclidean distances, $\Delta_p \equiv \|{\bf
s}_{i}-{\bf s}_{j}\|$, of all the vectors in $\mathfrak{B}_0$ i.e.,
\begin{equation}
\label{hata} \hat{a}_1^{d}=
\frac{1}{N_0} {\sum}_{p=1}^{N_0} \Delta_{p}^{d},
\end{equation}

This definition implies that $\hat{a}_1$ is a random variable that
depends on the sampling point configuration. In connection with the
consistency principle, the bandwidth $h_1$ is also a random
variable. However, since $\hat{a}_1$ represents an average over all
the near neighbor distances for all the points, its fluctuations are
not very significant. In particular, the coefficient of variation
declines with the number of sample points. To avoid cumbersome
notation we will not distinguish between $h_1^{m}$ and
$\hat{a}_1^{m}$, $m \in \mathbb{Z}^{+}$ and the respective
stochastic moments in the following theorems on the asymptotic
properties.

\begin{rem}{\rm
 Other estimators of the distance step
such as the median or the root mean square neighbor distances can be
used. However, the equation~(\ref{hata}) leads to consistent
convergence properties for the variance of the sample constraints,
regardless of the spatial dimension.}
\end{rem}

 The kernel bandwidth,
$\hat{h}_1 := h_1(\hat{a}_1),$ is then  given in view of
(\ref{eq:bandw1}) by

\begin{equation}
\label{eq:estim-h1} \hat{h}_1 =\hat{a}_1 \, B_2^{-1/2}.
\end{equation}
The above gives an explicit linear solution for the bandwidth in
terms of the step. In practical applications $p_n$ is a small
parameter, and (\ref{eq:estim-h1}) is sufficient. Alternatively,
(\ref{eq:alpha1}) can be solved numerically to obtain the bandwidth
in the pre-asymptotic case.

 \begin{lem}\label{lem-h1as}
 Let us assume that the sampling network densification conditions (\ref{H:band})
 and (\ref{H:dens}) hold, i.e., $|\Om_n| \propto n^{\delta}$, where
 $0<\delta <1$ and $h_{1} \propto  n^{-\gamma}$ for every realization
 of the sampling network. Then,
 if the gradient step is defined by the
 equation (\ref{hata}), the
 bandwidth exponent satisfies the inequality
 $0<\gamma \le (1-\delta)/d$.
 \end{lem}

{\bf Proof:} It holds that $\sum_{p=1}^{N_0} \Delta_p^d \geq
v_d\,|\Om_n|$ (where $v_d$ is a geometric constant that depends on
$d$). Since $N_0>n$, in light of equations (\ref{hata}) and
(\ref{eq:estim-h1}) it follows that $\hat{h}_1^d > v_{1,d}' \,
n^{-(1-\delta)}$ for any $n$, which implies that $d\,\gamma \leq
1-\delta.$


\begin{theo}{{\rm (Mean of the Sample Gradient Constraint - Differentiable Case.)}}
\label{theor:grad-dif} Assume that conditions
(\ref{H:loc})-(\ref{H:dens}) above are satisfied, and that
$F_{\la}$ is four times differentiable in a neighborhood of zero.
Then $\phiso$ is an asymptotically unbiased estimator of the
stochastic constraint $\phieo$. More specifically, the following
holds

\begin{equation}\label{eq:theor1}
    \E \left[\overline{\varphi_1}(h_1) -\phi_1(a_1) \right] =
 d \,\tau_2\,h_1^4 \left( B_4 - B_2^{2} \right)
 +O(h_1^2\,p_n)+o(h_1^4),
\end{equation}
where $\tau_2=F_{\la}^{(2)}(0)/2.$
\end{theo}

{\bf Proof:}
Since $\phiso=d\,\fbarg$ and $\phieo=d \,F_{\la}(a_1)$ we focus on
$\fbarg$ and $F_{\la}$ to avoid unnecessary clutter. The sample
function
 $\fbarg$, defined in (\ref{fbar}), is expressed as follows in light of
 equations~(\ref{Ksum}) and~(\ref{kernel-av}) :

\begin{equation}
\label{eq:KhS1-W} \fbarg= \frac{ {\sum}'_{i, j} W_{i,j} \,
      \left( \Xo_{i,j}\right)^2 }
     {  {2\sum}'_{i, j} W_{i,j}  }.
\end{equation}
$\E[\fbarg]$ involves an expectation over both the sampling point
distribution and the distribution of the field values. Hence, we can
write
$$\E[\fbarg]=\E \left\{ \E \left [\fbarg/\bfU_1,\ldots,\bfU_n\right
]\right\},$$ where the inner (conditional) expectation is over the
field values keeping the sampling locations fixed, whereas the outer
expectation is with respect to the sampling point distribution.

\vspace{6pt}

\textit{Calculation of  the Conditional Expectation $\E \left
[\fbarg/\bfU_1,\ldots,\bfU_n\right ].$} \vspace{4pt}

Since only the numerator of (\ref{eq:KhS1-W}) depends on the field
values, we obtain

\begin{equation}
\label{eq:Econd-fh1} \E \left [\fbarg/\bfU_1,\ldots,\bfU_n\right ]=
\frac{ \ka_{h_1} \Big\{ F_{\la}(\vth \bfU_{i,j})\Big\} } {\ka_{h_1}
\left\{ 1\right\}}.
\end{equation}

\vspace{6pt}

\textit{Leading-order calculation of $\ka_{h_1} \Big\{ F_{\la}(\vth
\bfU_{i,j})\Big\}.$} \vspace{4pt} \\
Using Lemma~(\ref{lem:ergod}) we obtain
\begin{eqnarray*}
 & &
 \E\left[ W_{i,j}\, \Big \{ F_{\la} \left(\vth \, \bfU_{i,j}
\right)\Big\}\right] =
 \int d\bfom\;K\big ( \|\bfom\| /p_n\big )
 \;F_{\la}\big (\vth \|\bfom\| \big) f_1(\bfom) \\
&=&  \int d \om \,\om^{d-1} \, K\big ( \om /p_n\big ) \, F_{\la}\big
(\vth \om \big)
 \int_{\mathbb{S}_d}d\bftheta\;
  J_d(\bftheta)f_1\left(\om \bfz \right) \\
   &=& p_n^{d}A_d\big [f_1(0 )+O(p_n)\big ]
   \int d u\,  u^{d-1} K\big ( u\big ) F_{\la}\big (h_1
u \big) .
 \end{eqnarray*}
 Since the kernel is compactly supported, it is possible to
 approximate $F_{\la}(h_1u) $ with a Taylor series expansion around
  zero, i.e.,
  \begin{equation}
  \label{eq:Fexdif}
 F_{\la}(h_1u)=\tau_2\,u^2h_1^2+\tau_4\,u^4h_1^4+o(h_1^4),
 \end{equation}
where $\tau_i=F_{\la}^{(i)}(0)/i!,$ $i=2,4$. Inserting the expansion
in the integral it follows that
\begin{equation*}
 \E\left[ W_{i,j}\, \Big \{ F_{\la} \left(\vth \, \bfU_{i,j}
\right)\Big\}\right]
   =p_n^{d}A_d\big [f_1(0 )+O(p_n)\big ] \, \left[ \tau_2\,h_1^2 m_{K,d+2}+
   \tau_4 \,h_1^4 m_{K,d+4} +o(h_1^4)\right].
 \end{equation*}
Finally, based on the above and Lemma~(\ref{lem:ergod}), it follows
that
\begin{equation*}
 \ka_{h_1} \Big\{ F_{\la}(\vth \bfU_{i,j})\Big\}=
n^2\,p_n^{d}A_d\big [f_1(0 )+O(p_n)\big ]\left[ \tau_2\,h_1^2
m_{K,d+2}+
   \tau_4 \,h_1^4 m_{K,d+4} +o(h_1^4)\right]
   \quad \mbox{a.s}.
\end{equation*}
   Using this equation in connection with (\ref{eq:ewij}) and (\ref{eq:KhS1-W})
   leads to

\begin{equation}
\label{eq:E-fh1} \E \left [\fbarg/\bfU_1,\ldots,\bfU_n\right ] =
\tau_2\,h_1^2 \, B_2 + \tau_4\,h_1^4 \, B_4
+O(h_1^2p_n)+o(h_1^4)\quad \mbox{a.s}.
\end{equation}
The respective expansion for $F_{\la}(a_1)$  is obtained using the
consistency principle, (\ref{eq:bandw1}), as well as
equations~(\ref{eq:Fexdif}) and (\ref{eq:estim-h1}), i.e.,

\begin{equation}
\label{eq:Fl-dif} F_{\la}(a_1)=\tau_2\,h_1^2 \, B_2 +\tau_4\,h_1^4
\, B_2^2 + O(h_1^2p_n)+o(h_1^4), \quad {\rm a.s.}.
\end{equation}
From the equations (\ref{eq:E-fh1}), (\ref{eq:Fl-dif}) and
Lemma~(\ref{lem:useful}), it follows that

\begin{equation}
\label{eq:fxh1}
 \E \left [\fbarg-F_{\la}(a_1)\right ]= \tau_4\,h_1^4 \, \left(
B_4 - B_2^2 \right) + O(h_1^2p_n) + o(h_1^4).
\end{equation} The asymptotic convergence then follows from the
densification effect, i.e., from $\gamma>0$. In light of the above,
$\fbarg$ is an asymptotically unbiased estimator of $F_{\la}(a_1),$
since the difference $\E \left [\fbarg\right ]-F_{\la}(a_1)$
converges to $0$ faster than each component as $h_1 \rightarrow 0.$

If we consider fluctuations in the bandwidth and the step,  $h_1$ on
the right hand-side of equation~(\ref{eq:fxh1}) should be replaced
by the respective mean value, and the corrections should also
include bandwidth fluctuations.

\begin{rem}{\rm The asymptotic decline of the bias as $h_{1}^{4}$
follows from the consistency principle and does not require the
specific choice of the step (\ref{hata}). The latter may only
influence the upper bound of the bandwidth exponent $\gamma.$}
\end{rem}

The following is also a direct consequence of
Theorem~(\ref{theor:grad-dif}) and Lemma~(\ref{alpha1-as}):

\begin{equation}
\label{eq:S1bias} \E[\overline{\mathcal{S}_1}(a_1)] -
\E[{S_1}(a_1)]= \tau_4\,h_1^2 \, \left(\frac{B_4}{B_2}-B_2\right) +
O(p_n) + o(h_1^2).
\end{equation}
Equation~(\ref{eq:S1bias}) shows that the generalized gradient
$\overline{\mathcal{S}_1}(a_1)$ is an asymptotically unbiased
estimator of the stochastic constraint $\E[{S_1}(a_1)]$.

\begin{theo}
\label{theor:grad-cont} {{\rm (Mean of the Sample ``Gradient''
Constraint - Continuous Case.)}} Assume that conditions
(\ref{H:loc})-(\ref{H:dens}) above are satisfied, and that $F_{\la}$
is continuous but non differentiable at zero. For $a'_1=h_1\,B_1$,
it follows that $\phiso$ is an
 is an asymptotically unbiased estimator of the
stochastic constraint $\phi_1(a'_1).$ More specifically:

\begin{equation}\label{eq:biasphi1-1}
    \E [\overline{\varphi_1}(h_1) -\phi_1(a'_1)  ]=  d \,\tau_2\,h_1^2
    \left(B_2- B_1^2\right)  + d \, \tau_3\,h_1^3 \, \left
(B_3 - B_1^{3} \right) + O(h_1\, p_n) + o(h_1^3).
\end{equation}

In addition, $\phiso$ is an
 asymptotically biased estimator of the
stochastic constraint $\phi_1(a_1),$ i.e.,

\begin{equation}\label{eq:biasphi1-2}
    \E [\overline{\varphi_1}(h_1)-\phi_1(a_1)   ] =
  d\, \tau_1\,h_1 \, \left(B_1-B_2^{1/2}\right ) +
  d \,\tau_3\,h_1^3 \, \left(B_3-B_2^{3/2}\right )
 + O(h_1p_n)+o(h_1^3).
\end{equation}
The mean relative error of $\E [\overline{\varphi_1}(h_1) ]$  is
given by

\begin{equation}
\label{eq:relerr1}
 \psi_{\epsilon,1} :=    \E \left[ \frac{\overline{\varphi_1}(h_1)
 -\phi_1(a_1)}{\phi_1(a_1)} \right]
  = \frac{B_{1,2}}{\sqrt{B_{2}}}
   +     h_1^2\,  \frac{\tau_3}{\tau_1}\frac{B_{3,2}}{\sqrt{B_{2}}}
   +O(p_n)+o(h_1^2),
\end{equation}
where $B_{1,2}=B_{1}-B_{2}^{1/2},$ and $B_{3,2}=B_{3}-B_{2}^{3/2}.$
\end{theo}

{\bf Proof:} The logic of the proof is the same as in
Theorem~(\ref{theor:grad-dif}), and therefore we only present the
main points.
 The derivatives of $F_{\la}$ do not exist at zero.
 However, if $F_{\la}$ admits at least third-order derivatives for any $h_1 u >0$,
 the Taylor series expansion of the semivariogram
 is expressed as

\begin{equation}
 \label{eq:Fexcon}
F_{\la}(h_1u)=\tau_1\, uh_1+\tau_2 \,u^2h_1^2+\tau_3\, u^3h_1^3 +
O(h_1\,p_n) + o(h_1^3), \end{equation} where
$\tau_i=F^{(i)}_{\la}(0^{+})/i!$.
 Then  we obtain
\begin{eqnarray*}
 \E\left[ W_{i,j}\, \Big \{ F_{\la} \left(\vth \, \bfU_{i,j}
\right)\Big\}\right] =p_n^{d}A_d\big [f_1(0 )+O(p_n)\big ]\left[
{\sum}_{j=1}^{3} \tau_j\,h^j m_{K,d+j}+o(h_1^3)\right],
 \end{eqnarray*}
 and in connection with (\ref{eq:ewij}) and (\ref{eq:KhS1-W}) it
 follows that

\begin{equation}
\label{eq:EE-fh1} \E \left [\fbarg/\bfU_1,\ldots,\bfU_n\right ]=
\sum_{j=1}^{3} \tau_j\,h^j \, B_j +o(h_1^3)\quad \mbox{a.s}.
\end{equation}
Based on (\ref{eq:Fexcon}), the semivariogram $F_{\la}(a'_1)$ is
expressed as
\begin{equation}
F_{\la}(a'_1)=\sum_{j=1}^{3}\tau_j\,h^j \, B_1^j
+O(h_1p_n)+o(h_1^3).
\end{equation}
Hence, we obtain

\begin{equation}
\E \left [\fbarg/\bfU_1,\ldots,\bfU_n \right ] - F_{\la}(a'_1)  =
\sum_{j=2,3} \tau_j\,h_1^j \left( B_j-B_1^j \right) +
O(h_1\,p_n)+o(h_1^3) \quad {\rm a.s}.
\end{equation}
The above proves
equation~(\ref{eq:biasphi1-1}).  The equation~(\ref{eq:biasphi1-2})
is proved in the same way, but the expansion~(\ref{eq:Fexcon}) is
replaced with an expansion around $a_1$, i.e.,

\begin{equation}
 \label{eq:Fexcon2}
F_{\la}(a_1)=\tau_1\, h_1 \, B_2^{1/2}+\tau_2\, h_1^{2} \, B_2 +
\tau_3\, h_1^{3} \, B_2^{3/2} + o(h_1^{3}).
\end{equation}
The above, in connection with~(\ref{eq:EE-fh1}),  leads to

\begin{equation}
 \label{eq:Fexcon3}
\E \left [\fbarg/\bfU_1,\ldots,\bfU_n \right ] - F_{\la}(a_1)  =
\sum_{j=1,3} \tau_j\,h_1^{j} \, \left(B_j - B_2^{j/2}\right )  +
O(h_1\,p_n)+o(h_1^3) \quad {\rm a.s}.
\end{equation}
 Finally, equation~(\ref{eq:relerr1}) for the mean relative error (relative bias),
follows from~(\ref{eq:Fexcon2}) and~(\ref{eq:Fexcon3}).

Based on equation~(\ref{eq:relerr1}) the relative bias depends on
the kernel function through the coefficients $B_{1,2}$ and
$B_{3,2}$. As $h_1 \rightarrow 0$, the relative bias converges
 to $B_{1,2}/B_2^{1/2}$.
\begin{lem}
\label{lem:relbias1} The asymptotic relative bias,
$\psi_{\epsilon,1},$ is a non-positive number.
\end{lem}
{\bf Proof:} $B_2$ is a positive number. By definition, $B_{1,2}=
B_1 - B_2^{1/2}$. Let us define the density function
$$f_K(s):=\frac{K(s) \, s^{d-1}}{\int_{0}^{R} K(s) \, s^{d-1}}, \quad
s \in [0,\,R].$$ In light of this definition and equation
(\ref{eq:mom-rat}), we obtain $B_m = \E_{K}[s^m],$ where $E_{K}$
denotes the expectation with respect to the density function $f_K.$
Then, $B_2 - B_1^{2} =\E_{K}\left[ \left( s- \E_{K}[s] \right)^{2}
\right] \ge 0$, and thus $B_{1,2} \le 0$ follows directly.

As a direct consequence of equations~(\ref{eq:estim-h1})
and~(\ref{eq:EE-fh1}), one obtains that
$\E[\overline{\mathcal{S}_1}(a_1)] \propto O(h_1^{-1})$. Hence, the
sample function $\overline{\mathcal{S}_1}(a_1)]$ is not well defined
at the asymptotic limit. Thus,  the ``gradient'' constraints in the
continuous but non-differentiable case refer to the sample function
$\overline{\varphi_1}(h_1)$ and its stochastic counterpart,
$\phi_1(a_1).$


 \begin{theo}{{\rm (Variance of the Sample ``Gradient'' Constraint.)}}
 \label{theor:grad-var} If the conditions (\ref{H:loc})-(\ref{H:dens})
 above are satisfied, $\overline{\varphi_1}(h_1)$
 is an asymptotically consistent
estimator of $\phi_1(a_1).$ In particular, the variance of $\phiso$
is given asymptotically by:

\begin{equation}
\label{varS1}
    \vrc \left[  \overline{\varphi_1}(h_1)\right]
  = O\left (\frac{1}{n^{2c_1 \gamma+\epsilon_1}}\right ), \quad
  \epsilon_1 = \min\{ \delta, 2-\delta-d\,\gamma\}.
\end{equation}
\end{theo}
{\bf Proof:}
The variance of $\fbarg$ is given by means of:

\begin{equation}
\label{eq:varS1-2} \vrc [\fbarg] = \E \Big [ \vrc
[\fbarg/\bfU_1,\ldots,\bfU_n] \Big]+ \vrc \Big [ \E
[\fbarg/\bfU_1,\ldots,\bfU_n] \Big]  .
\end{equation}

\noindent According to Eq.~(\ref{eq:E-fh1}) in the differentiable
case, and Eq.~(\ref{eq:EE-fh1}) in the non-differentiable case, the
second term on the right hand side of Eq.~(\ref{eq:varS1-2}) is
$o(h_1^{2}).$ Hence, we focus on the first term, which is expressed
as follows:

\begin{eqnarray}
\label{varS1-sum}  \vrc [\fbarg/\bfU_1,\ldots,\bfU_n] &= &
{\sum}'_{i, j}\,{\sum}'_{k, l} \frac{W_{i,j} \, W_{k,l} \,} {\left[
\ka_{h_1} ( 1)\right]^{2}}  \, \cvc \left \{
       (X^{*}_{i,j})^2 , \, (X^{*}_{k,l})^2  \right \} \nonumber \\
& =&  V_{1,1} + V_{1,2} + V_{1,3},
\end{eqnarray}
 where the functions $V_{1,1}, V_{1,2}, V_{1,3}$, in light of
 $\psi_{\la}$ defined in ~(\ref{psila}), are given by
\begin{eqnarray}
\label{eq:v11}
  V_{1,1}
  &=& 2\, \frac{
{\sum}'_{i,j} W_{i,j}^2  \,
  \psi_{\la}\left(0,0,s_{i,j},s_{i,j}\right) }
  { \left[ {\sum}'_{i, j} W_{i,j} \right]^{2} }
  \\
\label{eq:v12}
  V_{1,2}
 &=& 4  \, \frac{
{\sum}'_{i,j,k} W_{i,j} \, W_{k,i} \,
  \psi_{\la}\left(s_{i,k},s_{j,i},0,s_{j,k}\right) }
  {  \left[ \sum'_{i, j} W_{i,j}  \right]^{2} }
  \\
\label{eq:v13}
  V_{1,3}
  &=& \frac{
{\sum}'_{i,j,k,l} W_{i,j} \, W_{k,l} \,
  \psi_{\la}\left(s_{i,k},s_{j,l},s_{i,l},s_{j,k}\right)}
  {  \left[ {\sum}'_{i,j} W_{i,j}   \right]^{2} }.
\end{eqnarray}

\vspace{4pt} \textit{Leading-order calculation of the denominator.}

The quantities $V_{1,1}, V_{1,2}, V_{1,3}$ in equations
(\ref{eq:v11})-(\ref{eq:v13}) have a common denominator, the
asymptotic behavior of which follows from (\ref{eq:ewij}). More
precisely, the following is true:
\begin{equation}
\label{eq:v11-den} \left[ {\sum}'_{i,j} W_{i,j} \right]^{2} =n^4 \,
p_n^{2d} \, \left[ A_d \, f_1(0) \, m_{K,d}\right]^2
  +o\left(n^4p_n^{2d}\right) \quad \mbox{a.s.}
\end{equation}

\vspace{4pt} \textit{Leading-order calculation of $ V_{1,1}.$}

Denote the numerators of (\ref{eq:v11})-(\ref{eq:v13}) by $N^{({\rm
v})}_{1,j}, j=1,2,3$. Then, it follows from Lemma~(\ref{lem:ergod})
that $N^{({\rm v})}_{1,1} = 2 n^{2} \, \tilde{N}^{({\rm v})}_{1,1}
\,$ almost surely, where :
\begin{eqnarray*}
\tilde{N}^{({\rm v})}_{1,1}&:=&\E \left[ W_{i,j}^{2} \,
\psi_{\la}\Big(0,0,\vth \, U_{i,j}, \vth \, U_{i,j} \Big)\right]\\
 & =& \int d\bfom \;  K^2\left  ( \|\bfom\|/ p_n \right )
   \psi_{\la}\Big(0,0,\vth \|\bfom\|,\vth \|\bfom\|\Big) f_1(\bfom).
\end{eqnarray*}
We use the variable $ u=\om / p_n$, and a Taylor expansion of $f_1$
around zero. We evaluate the integral over $u$ with the mean value
theorem. Finally, we apply the first scaling property of
$\psi_{\la}$ in condition (\ref{H:psil-prop}), to obtain the
following
\begin{eqnarray*}
 \tilde{N}^{({\rm v})}_{1,1}
  &=&  p_{n}^{d}   \int du\; u^{d-1} \,
K^2(u) \, \psi_{\la}(0,0,uh_1,uh_1)\int_{\mathbb{S}_d}d\bftheta\;
  J_d(\bftheta)\Big[ f_1(0 )+o(1)\Big]
\nonumber \\
&=& p_{n}^{d}\,g_1(u^*)\left (\frac{ h_1}{\xi}\right )^{2c_1}
A_d\,f_1(0) \int_{0}^{R} du\; u^{d-1} \, K^2(u)
+o(p_{n}^{d}\,h_1^{2c_1}).
\end{eqnarray*}
Hence, it follows that $N^{({\rm v})}_{1,1}$ is given by
\begin{equation}
\label{eq:v11-num} N^{({\rm v})}_{1,1} = 2 \,
n^2\,p_{n}^{d}\,g_1(u^*)\left (\frac{ h_1}{\xi}\right )^{2c_1}
A_d\,f_1(0) \, m^{(2)}_{K,d} + o(n^2\,p_{n}^{d}\,h_1^{2c_1}) \quad
{\rm a.s.}
\end{equation}
Finally, from equations (\ref{eq:ewij}), (\ref{eq:v11-den}),
(\ref{eq:v11-num}) and based on Lemma~(\ref{lem:ergod}), it follows
that
\begin{equation}
\label{eq:v11-res} V_{1,1}=
  \frac{2\,g_1(u^*)}{A_d \, f_1(0) \, \xi^{2c_1}} \, \frac{ m^{(2)}_{K,d}}{(m_{K,d})^{2}}
  \, \left( \frac{\vth^d}{ n^2 \, h_1^{d-2c_1} } \right)
   +
  o\left (\frac{\vth^d}{ n^2 \, h_1^{d-2c_1} } \right ) \quad {\rm
  a.s.}
\end{equation}
Hence, the asymptotic  dependence of $V_{1,1}$ on $n$ becomes

\begin{equation}
\label{eq:v11-n} V_{1,1}= O(n^{\delta-2+\gamma d-2c_1\, \gamma}).
\end{equation}

\vspace{4pt} \textit{Leading-order calculation of $ V_{1,2}.$}

The numerator of $V_{1,2}$ is equal to $N^{({\rm v})}_{1,2} = 4
n^{3} \, \tilde{N}^{({\rm v})}_{1,2} $ almost surely, where:
\begin{eqnarray*}
\tilde{N}^{({\rm v})}_{1,2} &:=& \E \Big[ K\left( U_{i,j}/
p_n\right) K\left( U_{i,k}/ p_n\right ) \,
 \psi_{\la}\left(\vth U_{i,k},\vth U_{i,j},0, \vth U_{j,k}\right)
 \Big]\\
  &=& \iint d\bfom_1 d\bfom_2\; K\left (\frac{\|\bfom_1\|}{p_n}\right )
 K\left (\frac{\|\bfom_2\|}{p_n}\right )
    \psi_{\la}\Big(\vth \|\bfom_2\|,\vth \|\bfom_1\|,0,\vth \|\bfom_{1,2}\|\Big) \,
    f_2(\bfom_1,\bfom_2).
\end{eqnarray*}
Converting $ \bfom_1$ and $\bfom_2$ to spherical polar coordinates,
 using the perturbation parameter $p_n$ with the change of variables
 $u_1 =\om_1/p_n$, $u_2 =\om_2/p_n$  leads to:
\begin{eqnarray*}
 \tilde{N}^{({\rm v})}_{1,2}
&=& p_n^{2d}\int du_1\; u_1^{d-1}\; K\left(u_1\right)
 \int du_2\; u_2^{d-1}\;
  K\left( u_2\right ) \prod_{i=1,2} \int_{\mathbb{S}_d}d\bftheta_i
  J_d(\bftheta_i)
    \nonumber \\
  & &
    \psi_{\la} \left(\vth \om_2,\vth \om_1,0,\vth \big\|\om_1\,\bfz_1
    -\om_2\,\bfz_2 \big \|\right)
    f_2 \left(\om_1\,\bfz_1, \om_2\,\bfz_2 \right).
  \end{eqnarray*}
We evaluate the integrals over $\bfom_1$ and $\bfom_2$ using the
mean value theorem, defining
 $u^{*}_{1,2}:= \big\|u_1^*\,\bfz_1^* -u_2^*\,\bfz_2^* \big \|$.
  By applying
the second scaling property of condition (\ref{H:psil-prop}) for
$\psi_{\la}$, the following is obtained:
\begin{eqnarray*}
\tilde{N}^{({\rm v})}_{1,2} &= &
    \psi_{\la} \left(h_1 u_2,h_1 u_1,0,h_1 \big\|u_1\,\bfz_1-u_2\,\bfz_2 \big \|\right)
    f_2 \left(p_n\, u_1\,\bfz_1, p_n \,u_2\,\bfz_2 \right)
        \nonumber \\
    &=& g_2 \left( u_2^*, u_1^*,0,u^{*}_{1,2}\right)\left (\frac{h_1}{\xi} \right)^{2c_1}
    \,p_n^{2d}\Big (m_{K,d}\, A_d\Big)^2
    f_2 \left(0 \right)+o\left (p_n^{2d} \, h_1^{2c_1}\right ).
\end{eqnarray*}
Hence, the following expression is obtained for the numerator of
$V_{1,2}$:
\begin{eqnarray}
\label{eq:v12-num} N^{({\rm v})}_{1,2}
  & = & 4 \, n^3  \,p_n^{2d} \, m^{2}_{K,d}\, A_d^2 \,
    f_2 \left(0 \right) \, g_2
  \left( u_2^*, u_1^*,0,u^{*}_{1,2}\right)\left (\frac{h_1}{\xi}
    \right)^{2c_1}
     \nonumber \\
    & +  &  o\left (p_n^{2d} \,n^3\,h_1^{2c_1}\right ) \quad {\rm
    a.s.}
\end{eqnarray}
Finally, based on equations (\ref{eq:ewij}), (\ref{eq:v11-den}),
(\ref{eq:v12-num}) and Lemma~(\ref{lem:ergod}), the following
asymptotic expression is obtained for $V_{1,2}$
\begin{equation}
\label{eq:v12-res} V_{1,2}= \frac{4 g_2
  \left( u_2^*, u_1^*,0,u^{*}_{1,2}\right)
    f_2 \left(0 \right) }{\xi^{2c_1}\,f^{2}_1(0)} \,
    \left (\frac{h_1^{2c_1}}{n}\right )
    + o\left (\frac{h_1^{2c_1}}{n }\right ) \quad {\rm a.s.}
\end{equation}

Therefore, the asymptotic  dependence of $V_{1,2}$ on $n$ becomes

\begin{equation}
\label{eq:v12-n} V_{1,2}= O(n^{-1-2c_1\, \gamma}).
\end{equation}

\vspace{4pt} \textit{Leading-order calculation of $ V_{1,3}.$}

The numerator of $ V_{1,3}$, $N^{({\rm v})}_{1,3}$, includes a
summation over quartets of sampling points and  thus involves the
joint pdf of three independent distances $\Big (\bfU_{i,j},
\bfU_{k,l}, \bfU_{i,k} \Big )$. For reasons of brevity, we denote $
{\bf u}_{i,j}=\bfom_1$, $ {\bf u}_{k,l}=\bfom_2 $, and $ {\bf
u}_{i,k}=\bfom_3$; then ${\bf u}_{i,l}= \bfom_2 + \bfom_3,$ and
${\bf u}_{j,l}= \bfom_2 + \bfom_3 -\bfom_1$; also
$u_{i,l}=\|\bfom_2+\bfom_3\|$ and
$u_{j,l}=\|\bfom_2+\bfom_3-\bfom_1\|$.

 According to Lemma (\ref{lem:ergod}), $N^{({\rm v})}_{1,3}=n^{4}
\tilde{N}^{({\rm v})}_{1,3}$ almost surely, where $\tilde{N}^{({\rm
v})}_{1,3}$ is given by

\begin{eqnarray*}
\tilde{N}^{({\rm v})}_{1,3} & := & \E \Big[ K\left(
U_{i,j}/p_n\right) K\left(U_{k,l}/p_n\right) \,
  \psi_{\la}\left( \vth U_{i,k},\vth U_{j,l},\vth U_{i,l},\vth U_{j,k}\right)
    \Big] \\
  & =& \iiint d\bfom_1   d\bfom_2  d\bfom_3 \;
  K\left(\frac{\|\bfom_1\|}{p_n}\right )
  K\left (\frac{\|\bfom_2\|}{p_n}\right ) \,
  f_3\left(\bfom_1,\bfom_3,\bfom_2+\bfom_3\right)\\
  & &
  \psi_{\la}\left (\vth \om_3,\vth u_{j,l},
  \vth u_{i,l}, \vth \om_{1,3}\right ).
\end{eqnarray*}
Converting $ \bfom_1,$ $\bfom_2$ and $\bfom_3$ to spherical polar
coordinates, the following expression is obtained:
\begin{eqnarray*}
\tilde{N}^{({\rm v})}_{1,3}
  & = &  \prod_{i=1}^{3} \int d\om_i   \, \big(\om_1 \om_2 \om_3\big)^{d-1} \,
  K\left (\frac{\om_1}{p_n}\right )
  K\left (\frac{\om_2}{p_n}\right )\prod_{i=1}^{3} \int_{\mathbb{S}_d}
  d\bftheta_i \, J_d(\bftheta_i) \\
  & &
  \psi_{\la}\left (\vth \om_3,\vth u_{j,l},
  \vth u_{i,l}, \vth \om_{1,3}\right )
  f_3\left(\om_1 \bfz_1,\om_3 \bfz_3,\om_2 \bfz_2+\om_3 \bfz_3\right).
\end{eqnarray*}
Using the variable transformations $u_1=\om_1 / p_n,$ $u_2=\om_2 /
p_n $, $u_3=\vth \, \om_3,$  and the Taylor expansion of $f_3$
around $(0,0,0)$, $N^{({\rm v})}_{1,3}$ is transformed as follows
\begin{eqnarray*}
\tilde{N}^{({\rm v})}_{1,3}
  & = &\frac{ p_n^{2d}}{\vth^d}\int du_1  \int du_2 \int du_3 \big(u_1 u_2 u_3\big)^{d-1}
  K\left (u_1\right ) \, K\left (u_2\right )\prod_{i=1}^{3} \int_{\mathbb{S}_d}
  d\bftheta_i \, J_d(\bftheta_i)  \\  & &
  \psi_{\la}\left ( u_3, \big\|h_1u_2\bfz_2+u_3\bfz_3-h_1u_1\bfz_1\big\|,
   \big\|u_3\bfz_3+h_1u_2\bfz_2\big\|, \right. \nonumber \\
  & & \left. \quad \,
   \big\|u_3\bfz_3-h_1u_1\bfz_1\big\|\right ) \, \Big
   [f_3(0,0,0)+o(1)\Big].
\end{eqnarray*}
The integrals over $u_1$, $u_2$ and $\bftheta_i$ are evaluated using
the mean value theorem and the third scaling property of condition
(\ref{H:psil-prop}):
\begin{eqnarray*}
\tilde{N}^{({\rm v})}_{1,3}
  &=& \frac{ p_n^{2d}}{\vth^d}\left (\frac{h_1}{\xi}\right)^{2c_1} f_3({\bf 0})\,
  \Big(m_{K,d}\,A_d\Big )^2\int du_3\,  u_3^{d-1} \, g_{3}
\left(u_3,u_1^*,u_2^*,\bftheta_1^*,\bftheta_2^*,\bftheta_3^* \right)
\nonumber \\
& &   +o(p_n^{2d}\vth^{-d}\,h_1^{2c_1}).
\end{eqnarray*}
Finally, the following expression is obtained for $N^{({\rm
v})}_{1,3}$

\begin{eqnarray}
\label{eq:v13-num} N^{({\rm v})}_{1,3}  & = C^{*} \, n^{4} \, \left(
\frac{ p_{n}^{2}}{\vth} \right)^{d} \,
  \left (\frac{h_1}{\xi}\right)^{2c_1} \,f_3({\bf 0})\,
  \left(m_{K,d}\,A_d \right)^2
  +o(p_n^{2d}\,\vth^{-d}h_1^{2c_1}),
\end{eqnarray}
where $C^{*}=\int du_3\,  u_3^{d-1}\,
  g_{3}\left(u_3,u_1^*,u_2^*,\bftheta_1^*,\bftheta_2^*,\bftheta_3^*
\right)$ is a finite constant thanks to
assumption~(\ref{H:g3bounds}).

Hence, based on equations (\ref{eq:esswij}), (\ref{eq:v11-den}),
(\ref{eq:v13-num}) and Lemma~(\ref{lem:ergod}), the following
asymptotic expression is obtained for $V_{1,3}$
\begin{equation}
\label{eq:v13-res} V_{1,3}=  \frac{C^* \,
    f_3 ({\bf 0}) }{\xi^{2c_1}\,f^{2}_1(0)}
    \left (\frac{h_1^{2c_1}}{\vth^{d} }\right )
     +  o\left (\frac{h_1^{2c_1}}{\vth^{d}}\right )
    \quad {\rm a.s.}
\end{equation}

Hence, the asymptotic  dependence of $V_{1,3}$ on $n$ becomes

\begin{equation}
\label{eq:v13-n} V_{1,3}= O(n^{-\delta - 2c_1\, \gamma}).
\end{equation}

%

\vspace{4pt} \textit{Variance Convergence Rate}. \vspace{4pt}

 Based on equations (\ref{eq:v11-n}), (\ref{eq:v12-n}) and
(\ref{eq:v13-n}), the convergence of $V_{1,3}$ is slower than  that
of $V_{1,2}$ since $\delta<1.$ The convergence of $V_{1,1}$ is
faster than that of $V_{1,3}$ if $\gamma \, d < 2(1-\delta)$. If
this condition holds, then $V_{1,3}$ is the rate-limiting term. In
light of Lemma~(\ref{lem-h1as}) this inequality is satisfied for the
bandwidth defined by~(\ref{hata}).

\begin{rem} {\rm The rate of convergence of the gradient estimator's
variance is the same for the differentiable and non-differentiable
cases. The three terms, i.e., $V_{1,1}, V_{1,2}, V_{1,3},$ possess
distinct convergence rates. These terms correspond to sample
functions that involve doublets, triplets and quartets of
non-identical sampling points. Using the step estimate (\ref{hata})
and the consistency principle, the slowest convergence rate
(asymptotically dominant term) is due to the term that involves
quartets of non-identical points. On intuitive grounds, we would
expect the same behavior to hold for different step estimates.}
\end{rem}

\subsection{Generalized Curvature Estimation}
\label{ssec:gen-curv} In this section we prove a relation between
the ``curvature'' step and the kernel bandwidth $h_2$, and we
propose an estimate for the step. We then investigate the asymptotic
properties of the mean and variance of the generalized curvature
estimator. In the process, we also show that to a first
approximation $\mu_1=\mu_2=1$ and we calculate the asymptotic
dependence of the leading corrections.

In the proofs of asymptotic dependence, we will use the following
modification of Lemma~(\ref{lem:ergod}).
\begin{equation}
\label{eq:sumQ}
  {\rm Pr}\Big( \lim_{n\rightarrow \infty} \ka_{rh_2}\{A_{i,j}\}
 = n^{2} \, \E\left[Q_{i,j(r)} \, A_{i,j}\right] \Big) = 1.
\end{equation}

\begin{lem}\label{alpha2-as}
The following relation holds between the bandwidth $h_2$ and the
curvature step $a_2$:

\begin{equation*}\label{eq:bandw2}
 a_2^{4}=h_2^{4} \, B_4 +O(h_2q_n)
 \quad {\rm a.s.}.
 \end{equation*}
 \end{lem}
{\bf Proof:}
 The proof is along the lines of Lemma~(\ref{alpha1-as}).
 According to the definition~(\ref{alpha2}) and the kernel-average equation (\ref{kernel-av}),
 the step $a_2$ is defined by:
\begin{equation}
\label{alpha2p1} a_2^4 :=\langle s^{4}_{i,j} \rangle_{h_2} =
\frac{\ka_{h_2} \left\{ s^{4}_{i,j} \right\}} {\ka_{h_2}
\left\{1\right\} }.
\end{equation}

\textit{Leading-order calculation of $\E[Q_{i,j (r)}].$}
\vspace{4pt}

\begin{eqnarray}
\label{eq:EQij} \E[Q_{i,j (r)}]&=&\int d\bfom \;
  K\big (\vth\|\bfom \|/r\,h_2\big) \, f_1(\bfom)
\nonumber \\
&=& \int d\om \;\om^{d-1} K\left( \om /r\,q_n\right)
 \int_{\mathbb{S}_d}d\bftheta\;
  J_d(\bftheta)f_1\left(\om \bfz \right)
  \\
  &=& r^d \, q_n^{d} \, A_d \, f_1(0) \, m_{K,d}
   + O\left(q_n^{d+1}\right).   \nonumber
\end{eqnarray}
Hence, we obtain
\begin{equation}
\label{esswij1} \ka_{rh_2} \left\{ 1 \right\}= n^2 \,r^d \, q_n^{d}
\, A_d \, f_1(0) \, m_{K,d}
  + O\left(n^2 \, q_n^{d+1}\right) \quad {\rm a.s.}.
\end{equation}

The term $\ka_{h_2} \left\{ s^{4}_{i,j} \right\}$ is a special case
of  $\ka_{rh_2} \left\{ s^{m}_{i,j} \right\}$, which we evaluate
below. \\

 \textit{Leading-order calculation of $\E\left[Q_{i,j
(r)}\, s_{i,j}^m\right].$} \vspace{4pt}
\begin{eqnarray*}
 \E\left[Q_{i,j (r)}\, s_{i,j}^m \right]  & = &
   \vth^m \, \int d\bfom \;\|\bfom \|^m \,
  K\left( \|\bfom \|/rq_n\right) \,  f_1(\bfom)
\nonumber \\
  & = & \vth^m \,r^{d+m} \, q_n^{d+m} \, A_d \, f_1(0) \, m_{K,d+m}
  + O\left(\vth^m \, q_n^{d+5}\right).
\end{eqnarray*}
Hence, it follows that
\begin{equation}
\label{esswij4} \ka_{rh_2} \left\{ s_{i,j}^m \right\}= n^2 \, \vth^m
\,r^{d+m} \, q_n^{d+m} \, A_d \, f_1(0) \, m_{K,d+m}
  + O\left(n^2 \, q_n^{d+1+m}\right) \quad {\rm a.s.}.
\end{equation}
From the  Eqs.~(\ref{esswij1}) and (\ref{esswij4}) it follows that
\begin{equation}
\label{esswij44}
 \langle s^{m}_{i,j} \rangle_{rh_{2}} \equiv \frac{\ka_{rh_2} \left\{ s_{i,j}^m \right\} }
{\ka_{rh_2} \left\{ 1 \right\} }= r^m\, h_2^m \, B_m  + O(h_2^m \,
q_n) \quad \mbox{a.s.}.
\end{equation} The asymptotic behavior of $a_2$ is obtained from (\ref{esswij44})
for $r=1$ and $m=4$:

\begin{equation}
\label{a2hat} a_2^4=h_2^4 \, B_4 + O(h_2^4 \, q_n) \quad {\rm a.s.}.
\end{equation}

The coefficients $\mu_1(h_1)$ and $\mu_2(h_2)$ appear in the
definition of the generalized curvature constraint.  We calculate
the asymptotic dependence of these coefficients.
\begin{lem}\label{lem:mu1mu2-as}
 The coefficients $\mu_1(h_1)$ and $\mu_2(h_2)$ are given
 asymptotically by:
 $$ \mu_j(h_j)=1+o(1)
 \quad {\rm a.s.,} \quad \mbox{for $j=1,2.$}$$
 \end{lem}

 {\bf Proof:}
Based on the equations (\ref{muhat1}) and (\ref{muhat2}) the
coefficients involve the averages $\langle
s^{2}_{i,j}\rangle_{rh_2}$ and $\langle s^{4}_{i,j}\rangle_{rh_2},$
where $r=1,2,\sqrt{2}$. Both averages are given by
equation~(\ref{esswij44}).  The lemma is proved following
straightforward but tedious algebraic manipulations.

For the curvature step we will use the same expression as for the
gradient step, i.e., $\hat{a}_2=\hat{a}_1$, given by
equation~(\ref{hata}).
 The kernel bandwidth, $\hat{h}_2 :=
h_2(\hat{a}_2),$ is then given in view of (\ref{a2hat}) as follows:


\begin{equation}\label{eq:bandw3}
 \hat{h}_2=\hat{a}_2 \, B_4^{-1/4}.
 \end{equation}

 \begin{lem}\label{lem-h2as}
  Let us assume that $h_{2} \propto  n^{-\nu}$.
 Then, if the curvature step is defined by the
 equation (\ref{hata}), the
 bandwidth exponent satisfies the inequality
 $0<\nu \,d \le 1-\delta$.
 \end{lem}
{\bf Proof:} The proof is completely analogous to the proof of
Lemma~(\ref{lem-h1as}) if $\gamma$ is replaced by $\nu$.

\begin{theo} {{\rm (Mean of the Sample Curvature Constraint - Differentiable Case.)}}
\label{theor:curv-mean}  Assume that  hypotheses
(\ref{H:loc})-(\ref{H:dens}) above are satisfied, and that $F_{\la}$
admits five derivatives in a neighborhood of zero. Then
$\overline{\varphi_2}(h_2)$  is an asymptotically unbiased estimator
of $\phi_2(a_2).$ More specifically, the following holds

\begin{equation}
\label{eq:biasphi2dif} \E \left[ \overline{\varphi_2}(h_2) -
\phi_2(a_2) \right] = -24d(d+4)\, \tau_6\,h_2^6\, \left( B_6
-B_4^{3/2} \right) + O(h_2^2 \, q_n)+o(h_2^6),
\end{equation}
where $\tau_6=$ and $g_k=m_{K,d+k}/m_{K,d}$.
\end{theo}

{\bf Proof:}

Based on (\ref{barS2irr}),  $\phi_2(a_2)$ is expressed in terms of
$\fbarcr$ as follows:
\begin{equation}
\label{eq:varphi2} \overline{\varphi_2}(h_2)=\frac{1}{2} \left[
\ctwo \mu_1(h_2) \,
 \fbarc-\cthree \mu_2(h_2)
 \fbarcc-\cone \fbarccc \right].
\end{equation}
The sample function $\fbarcr$ is defined in terms of (\ref{fbar}),
and it is expressed in light of (\ref{Ksum2}) as follows:
\begin{equation}
\label{eq:KhS2-W} \fbarcr= \frac{ {\sum}'_{i, j} Q_{i,j (r)} \,
      \left( \Xo_{i,j}\right)^2 }
     {  {2\sum}'_{i, j} Q_{i,j}  }.
\end{equation}
 Hence, $\E  [\overline{\varphi_2}(h_2)]$ is expressed in terms of
 $\E[\fbarcr]$. As in Theorem~(\ref{theor:grad-dif}), the
 ensemble average implies
 $\E[\fbarcr]=\E \left\{ \E \left [\fbarcr/\bfU_1,\ldots,\bfU_n\right
]\right\}.$

\vspace{6pt}

\textit{Calculation of  the Conditional Expectation $\E \left
[\fbarcr/\bfU_1,\ldots,\bfU_n\right ].$} \vspace{4pt}

Only the numerator of (\ref{eq:KhS2-W}) depends on the field values,
i.e.,
\begin{equation}
\label{eq:Econd-fh2} \E \left [\fbarcr/\bfU_1,\ldots,\bfU_n\right ]=
\frac{ \ka_{rh_2} \Big\{ F_{\la}(\vth \bfU_{i,j})\Big\} }
{\ka_{rh_2} \left\{ 1\right\}}.
\end{equation}

\vspace{6pt}

\textit{Leading-order calculation of $\ka_{rh_2} \Big\{ F_{\la}(\vth
\bfU_{i,j})\Big\}.$} \vspace{4pt}
\begin{eqnarray*}
 & &
 \E\left[ Q_{i,j (r)}\, \Big \{ F_{\la} \left(\vth \, \bfU_{i,j}
\right)\Big\}\right] =
 \int d\bfom\;K\left( \|\bfom\| /(rq_n)\right)
 \;F_{\la}\left(\vth \|\bfom\| \right) f_1(\bfom) \\
&=&  \int d \om \,\om^{d-1} \, K\left( \om /(rq_n)\right) \,
F_{\la}\left(\vth \om \right)
 \int_{\mathbb{S}_d}d\bftheta\;
  J_d(\bftheta)f_1\left(\om \bfz \right) \\
   &=& r^d\,q_n^{d}A_d\big [f_1(0 )+O(p_n)\big ]
   \int d u\,  u^{d-1} K(u) F_{\la}(r\,h_2u) .
 \end{eqnarray*}
 The sixth-order Taylor series expansion of
 $F_{\la}(r\,h_2u) $ around zero yields
 \begin{equation}
 \label{eq:Fexdif6}
F_{\la}(r\,h_2u)=\tau_2 \,r^2 u^2 h_2^2 + \tau_4\,r^4 u^4 h_2^4 +
 \tau_6\, r^6 h_2^6 + o(h_2^6).
 \end{equation}
Inserting the expansion in the kernel integral, it follows that
\begin{equation*}
\E\left[ Q_{i,j (r)}\, \Big \{ F_{\la} \left(\vth \, \bfU_{i,j}
\right)\Big\}\right]
    =  q_n^{d} \, A_d \, \left[f_1(0 ) + O(p_n)\right] \,
   \left[ {\sum}_{i=2,4,6} m_{K,d+i} \, \tau_i\,(r \, h_2)^i +o(h_2^6) \right].
 \end{equation*}
The above in connection with (\ref{esswij1}) for the kernel average
$\ka_{rh_2} \left\{ 1\right\}$ lead to:
\begin{equation}
\label{eq:E-fh2} \E \left [\fbarcr/\bfU_1,\ldots,\bfU_n\right ] =
\sum_{i=2,4,6} g_i \, \tau_i\,(r \, h_2)^i
 + O(h_2^2 \, q_n)+o(h_2^6) \quad {\rm a.s}.
\end{equation}
From (\ref{eq:E-fh2}) and (\ref{eq:varphi2}) it follows that the
$O(h_2^{2})$ term vanishes if the coefficients $\mu_1(h_1),
\mu_2(h_2)$ are defined as in equations (\ref{muhat1}) and
(\ref{muhat2}). Finally, we obtain
$$\E \left [\overline{\varphi_2}(h_2)/\bfU_1,\ldots,\bfU_n\right]=
-\czero \, B_4 \, \tau_4\,h_2^4 \, -24d(d+4)\, B_6 \, \tau_6\,h_2^6
+ O(h_2^2 \,q_n)+o(h_2^6) \quad {\rm a.s}.$$

Using the definition of $ \phi_2(a_2)$, equation (\ref{sto3}), the
expansion (\ref{eq:Fexdif6}), and the step - bandwidth relation,
(\ref{eq:bandw3}), a series expansion is obtained for $ \phi_2(a_2)$
$$ \phi_2(a_2)=-\czero \tau_4\,h_2^4 \, B_4 - 24d(d+4)\tau_6\,h_2^6\, B_4^{3/2}
+O(h_2^2q_n)+o(h_2^6), \quad {\rm a.s.}.
$$
The two preceding expansions allow calculating the bias for the
curvature constraint by subtracting the terms on the respective
sides.  The  proof is completed by applying Lemma \ref{lem:useful}
to obtain equation (\ref{eq:biasphi2dif}).

\begin{theo}
\label{theor:curv-cont} {{\rm (Mean of the Sample ``Curvature''
Constraint - Continuous Case.)}} Assume that conditions
(\ref{H:loc})-(\ref{H:dens}) above are satisfied, and that $F_{\la}$
is continuous but non differentiable at zero. For $a'_2=h_2\,B_1$,
it follows that $\phist$ is an asymptotically unbiased estimator of
the stochastic constraint $\phi_2(a'_2).$ More specifically, if
$c^{(4)}_{d} = \left[ \ctwo- \sqrt{2}\,\cthree-2\cone \right]$, and
$c^{(5)}_{d} = \left[ \ctwo-2\sqrt{2}\,\cthree-8\cone \right]$ then:
\begin{equation}\label{eq:biasphi2-1}
    \E \left[\overline{\varphi_2}(h_2)  - \phi_2(a'_2) \right] =
    c^{(5)}_{d} \, \tau_3\,h_2^3\, \left( B_3- B_1^{3} \right) +o(h_2^3)+O(h_2\, q_n).
\end{equation}

In addition, $\phist$ is an
 asymptotically biased estimator of the
stochastic constraint $\phi_2(a_2),$ i.e.,
\begin{eqnarray}\label{eq:biasphi2-2}
    \E \left[\overline{\varphi_2}(h_2) -\phi_2(a_2)  \right] & = &
  c^{(4)}_{d}\,
 \tau_1\,h_2 \,\left( B_1 -B_4^{1/4}\right )
  +     c^{(5)}_{d} \,\tau_3\,h_2^3\, \left( B_3 - B_4^{3/4} \right)
  \nonumber \\
  & + &  o(h_2^3)+O(h_2\,q_n).
\end{eqnarray}
The mean relative error of $\E [\overline{\varphi_2}(h_2) ]$  is
given by

\begin{equation}
\label{eq:relerr2}
 \psi_{\epsilon,2} :=    \E\left[ \frac{ \overline{\varphi_2}(h_2)
 -\phi_2(a_2)}{\phi_2(a_2)} \right]
  = \frac{B_{1,4}}{B_{4}^{1/4}}
   +     h_2^2\,  \left( \frac{\tau_3}{\tau_1} \right) \frac{B_{3,4}}{B_{4}^{1/4}}
   +O(p_n)+o(h_2^2),
\end{equation}
where  $B_{1,4}=B_{1}-B_{4}^{1/4},$ and $B_{3,4}=B_{3}-B_{4}^{3/4}.$
\end{theo}

{\bf Proof:}


First we calculate $\E \left
[\overline{\varphi_2}(h_2)/\bfU_1,\ldots,\bfU_n\right]$. This
requires calculation of $\E \left
[\fbarcr/\bfU_1,\ldots,\bfU_n\right ]$.  The latter is given in
equation (\ref{eq:Econd-fh2}). On the right hand side of that
equation, the denominator, $\ka_{rh_2}\{1\}$, is given by equation
(\ref{esswij1}).  The numerator, $\ka_{rh_2}\{F_{\la} \left(\vth \,
\bfU_{i,j} \right) \}$, converges to $n^2 \,\E\left[ Q_{i,j (r)}\,
\Big \{ F_{\la} \left(\vth \, \bfU_{i,j} \right)\Big\}\right]$
according to (\ref{eq:sumQ}).

 For $h_2 u
>0$ the Taylor expansion of $F_{\la}$ is given by
\begin{equation}
\label{eq:Fexcon-curv} F_{\la}(r\,h_2u)=\tau_1\, ruh_2+\tau_2
\,r^2u^2h_2^2+\tau_3\, r^3u^3h_2^3+o(h_2^3). \end{equation} where
$\tau_{i}= F^{(i)}_{\la}(0^{+})$.  Then, we obtain by the standard
procedure

 \begin{eqnarray*}
 \E\left[ Q_{i,j (r)}\, \Big \{ F_{\la} \left(\vth \, \bfU_{i,j}
\right)\Big\}\right] =q_n^{d}\,A_d \big [f_1(0 )+O(q_n)\big ]\left[
{\sum}_{j=1}^{3} \tau_j\,r^{j} \, h_2^j \,
m_{K,d+j}+o(h_2^3)\right].
 \end{eqnarray*}
Based on the above, it follows that

\begin{equation}
\label{eq:EE-fh2} \E \left [\fbarcr/\bfU_1,\ldots,\bfU_n\right ] =
\sum_{j=1,2,3} \tau_j\,r^j\,h_2^{j} \,B_j + \ O(h_2 \, q_n) +
o(h_2^3)\quad \mbox{a.s}.
\end{equation}
Finally, using Lemma~(\ref{lem:useful}) and equation
(\ref{eq:varphi2}), we obtain the following
\begin{eqnarray*}\E \left
[\overline{\varphi_2}(h_2)/\bfU_1,\ldots,\bfU_n\right]&=& c^{(4)}_d
\, \tau_1\,h_2 \, B_1  + c^{(5)}_d \, \tau_3\,h_2^3\, B_3
+o(h_2^3)+O(h_2q_n) \quad \mbox{a.s.},
\end{eqnarray*}
where the term $O(h_2^{2})$ vanishes due to cancelation of the
coefficients.  Based on (\ref{sto3}) and the expansion
(\ref{eq:Fexcon-curv}), the stochastic term is expressed as
$\phi_2(a'_2=B_1\,h_2)=c^{(4)}_d \, \tau_1\,h_2 \, B_1  + c^{(5)}_d
\, \tau_3\,h_2^3\, B_1^{3} +o(h_2^3).$  This expansion in connection
with the one above for $\E \left
[\overline{\varphi_2}(h_2)/\bfU_1,\ldots,\bfU_n\right]$ leads to
\begin{eqnarray*}
\E \left
[\overline{\varphi_2}(h_2)/\bfU_1,\ldots,\bfU_n\right]-\phi_2(a'_2)
& = & c^{(5)}_d \tau_3\,h_2^3\, \left( B_3 - B_1^{3}\right )
+o(h_2^3)+O(h_2q_n) \quad \mbox{a.s}.
\end{eqnarray*}
The proof of equation (\ref{eq:biasphi2-1}) is completed by applying
Lemma \ref{lem:useful} to the above result. The estimator is
asymptotically unbiased since the bias converges to $0$ faster than
either the sample or the stochastic constraints.

Equations (\ref{eq:biasphi2-2}) and (\ref{eq:relerr2}) follow along
the same lines. The main difference is that the stochastic
constraint now becomes $\phi_2(a_2)$, where $a_2= (B_4^{1/4})\,h_2$
according to (\ref{eq:bandw3}).


\begin{lem}
\label{lem:relbias2} The asymptotic relative bias,
$\psi_{\epsilon,2},$ is non-positive.
\end{lem}
{\bf Proof:} As $h_2 \rightarrow 0$, the relative bias converges
 to $B_{1,4}/B_4^{1/4}$.  $B_4$ is a positive number. By definition, $B_{1,4}= B_1 -
B_4^{1/4}$. Using the density function defined in
Lemma~(\ref{lem:relbias1}), we can write $B_4 - B_1^{4} \ge
\E_{K}\left[ \left\{ s^{2}- \E^{2}_{K}[s] \right\}^{2} \right] \ge
0$, from which it follows that $B_{1,4} \le 0$.

 \begin{theo}{{\rm (Variance of the Sample ``Curvature'' Constraint.)}}
 \label{theor:curv-var} If the hypotheses (\ref{H:loc})-(\ref{H:dens}) above are
satisfied, then $\overline{\varphi_2}(h_2)$   is an asymptotically
consistent estimator $\phi_{2}(a_2).$ More specifically, the
following holds:

\begin{equation}
\label{varS2}
    \vrc \left[  \overline{\varphi_2}(h_2)\right]
  = O\left (\frac{1}{n^{2c_1 \nu+\epsilon_2}}\right ), \quad
  \epsilon_2 = \min\{ \delta, 2-\delta-d\,\nu \}.
\end{equation}

\end{theo}

{\bf Proof:} The proof is based on the same approach as in
Theorem~\ref{theor:grad-var}. The calculations are more extended due
to the cross-products between the sample functions $\fbarc$,
$\fbarcc$ and $\fbarccc$. However, in this case we also obtain terms
containing doublets, triplets and quartets of sampling points. Since
the complications are of a trivial nature, the lengthy calculations
will be omitted here. Using for the curvature step equation
(\ref{hata}), the quartet term dominates the convergence. This term
leads to the slow asymptotic decline of the variance as
$O(h_2^{2c_1}/\vth^d)$ or equivalently as $O(n^{-2c_1\,\nu -\delta}
)$.


\subsection{Calculation of Asymptotic Bias} \label{ssec:asym-bias}

\begin{table}
\begin{tabular}{ccccc}
    & Triangular & Quadratic & Gaussian & Tricube \\ \hline \hline
  $B_{1,2}$ & $-0.0472$ & $-0.0440$ & $-0.1138$ & $-0.0390$ \\ 
  $B_2$ &
  $1/5$ & $1/3 $& $1$ & $22/91$ \\ 
  $\psi_{\epsilon,1}$ & $-0.1056$ & $-0.0762$ & $-0.1138$ & $-0.0793$ \\ 
   \hline
  $B_{1,4}$ & $-0.1170$ & $-0.1056$ & $-0.3030$& $-0.0959$ \\ 
  $B_4 $ &
  $1/14$ & $1/6$ & $2$ & $22/243$  \\ 
  $\psi_{\epsilon,2}$ & $-0.2263$ & $-0.1653 $& $-0.2548$& $0.1748$ \\
  \hline \\
\end{tabular}
\caption{Calculations of $B_{1,2}$, $B_{1,4}$, $B_2$, $B_4 $, and
the relative bias of the ``gradient'' and ``curvature'' constraint
estimators using different kernel functions.} \label{Table1}
\end{table}
The asymptotic relative bias of the ``gradient'' and ``curvature''
constraint estimators obtained for different types of kernel
functions, according to equations~(\ref{eq:relerr1}) and
(\ref{eq:relerr2}), is shown in Table~(\ref{Table1}). In particular,
we include the {\it Gaussian kernel}, $K(s)=\exp(-s^2)$, the {\it
triangular kernel}, $K(s)=(1-\|s\|) \, \bigone_{ 0\leq \|\bfs\| \leq
1}$, the {\it quadratic kernel}, $K(s)=(1-\|s\|^{2}) \, \bigone_{
0\leq \|\bfs\| \leq 1}$, and the {\it tricube kernel},
$K(s)=(1-\|s\|^{3})^{3} \, \bigone_{ 0\leq \|\bfs\| \leq 1}$.  The
Gaussian kernel is not compactly supported, but it decays to zero
very fast. The quadratic kernel gives the lowest relative bias,
followed by the tricube kernel.

\section{SSRF Model Parameter Inference}
\label{sec:inference} Model parameter inference is based on the
procedure introduced in \cite{dth03}, which is expanded herein. The
main idea is to estimate the SSRF parameters, ${\bm \theta}$ by
matching sample constraints with their stochastic counterparts. We
use the sample constraints $ \overline{\mathcal{S}_0},$ $
\overline{\mathcal{S}_1}(\hat{a}_1)$ and $
\overline{\mathcal{S}_2}(\hat{a}_2)$, given by
equations~(\ref{barS0}), (\ref{barS1irr}), (\ref{barS2irr})
respectively,
 as well as the stochastic constraints $\E[S_0],$
 $\E[S_1(\hat{a}_1)]$ and
$\E[S_2(\hat{a}_2)], $ given by equations~(\ref{stoch-0}),
(\ref{semi-der2}), (\ref{semi-der4}) respectively. We also impose
the normalization constraint (\ref{s0prime}).

Determining ${\bm \theta}$ is then expressed as an optimization
problem that aims at minimizing the deviation between the stochastic
moments and their estimators; the latter include the sample-based
variance, gradient, and curvature constraints, and  $1$ for the
normalization constraint.  We introduce the
 metric $\Phi\left({\bf X}^{\ast};{\bm \theta}' \right)$
 to measure the  distance  between the sample and ensemble constraints:

\begin{eqnarray}\label{functional2}
    \Phi\left({\bf X}^{\ast};{\bm \theta}' \right) & := &
    \left\{1-\left (S_0'\right )^{1/\beta}\right\}^2+
    \left\{ 1-\left (
    \frac{\overline{\mathcal{S}_0}}{\overline{\mathcal{S}_1}(\hat{a}_1)}
    \frac{\E[S_1(\hat{a}_1)]}{\E[S_0]}\right )^{1/\beta} \right\}^2 \nonumber \\
    & + & \left\{ 1-\left(\frac{\overline{\mathcal{S}_1}(\hat{a}_1)}
    {\overline{\mathcal{S}_2}(\hat{a}_2)}
    \frac{\E[S_2(\hat{a}_2)]}{\E[S_1(\hat{a}_1)]}\right)^{1/\beta} \right\}^2\end{eqnarray}
In equation~(\ref{functional2}), $\E[S_0],$ $\E[S_1(\hat{a}_1)]$ and
$\E[S_2(\hat{a}_2)]$ are the values of the constraints obtained for
the ``current'' values of the Spartan parameters $\eta_1,$ $\xi$ and
$k_{\rm max}$ and for $\hat{a}_1=\hat{a}_2$ as given by
(\ref{hata}). The simplex search method of Nelder and Mead
\cite{simplex} is used for the optimization. The initial parameter
vector ${\bm \theta^{(0)}}$ is updated at every optimization step.
For $\eta_1$ the value $\eta_1^{(0)}=1$ is arbitrarily chosen. The
initial value $\xi^{(0)}$ of the characteristic length is estimated
from the data. The initial estimate of the characteristic length is
given by $\xi^{(0)}=\left[
\overline{\mathcal{S}_1}(\hat{a}_1)/\overline{\mathcal{S}_2}(\hat{a}_2)\right]^{1/2}.
$ The frequency cutoff $\km$ is chosen according $\km^{(0)}=2\pi
\hat{a}_1^{-1}$. The final vector, $\hat{\bm \theta}$, to which the
optimization converges gives the optimal parameters of the SSRF
model.

Note that the functional $\Phi\left({\bf X}^{\ast};{\bm \theta}'
\right)$ is independent of $\eta_0$, which can be set equal to $1$
during the optimization. The optimal $\hat{\eta}_{0}$ is obtained
using the condition of the variance independence from $\km \xi$ and
$\eta_1$, i.e., equation (\ref{eq:eta0}), which leads to the
following:
\begin{equation}\label{eta0optimal}
   \hat{\eta}_{0}=2^d\,\pi^{d/2} \Gamma(d/2)\, \hat{\sigma}_{\rmx}^2.
\end{equation}
In light of (\ref{stoch-0}) and (\ref{s0prime}),
equation~(\ref{eta0optimal}) guarantees that the model variance
matches the sample variance.

Some comments are in order regarding the definition of the distance
metric (\ref{functional2}). The functional is of the general form
$\Phi=\sum_{i=1}^{3} (1-z_i^{1/\beta})^{2}$, where
$$ z_1=S'_{0}, \quad z_2=\frac{\overline{\mathcal{S}_0}\,
\E[S_1(\hat{a}_1)]}{\overline{\mathcal{S}_1}(\hat{a}_1) \, \E[S_0]},
\quad z_3=\frac{\overline{\mathcal{S}_1}(\hat{a}_1) \,
\E[S_2(\hat{a}_2)]}
{\overline{\mathcal{S}_2}(\hat{a}_2)\,\E[S_1(\hat{a}_1)]}.$$ The
number of terms (squares) in $\Phi$ is equal to the number of
variables. The $z_i$ are functions that involve specific sample and
ensemble constraints. The definitions of $z_2$, $z_3$ are motivated
by the goals of (i) eliminating the dependence on $\eta_0$ (since
the latter is an overall scaling factor) (ii) defining dimensionless
variables so that the optimization does not depend on the units used
and (iii) forming combinations of constraints of similar magnitude
so that they contribute on an equal footing in the optimization.

Straightforward constraint differences, i.e.,
$\overline{\mathcal{S}_i} - \E[S_i]$ are neither dimensionless nor
of similar magnitude. Using ratios
$\overline{\mathcal{S}_i}/\E[S_i]$ yields dimensionless ratios of
similar magnitude, but it preserves the $\eta_0$ dependence. The
proposed combinations for $z_i$ for $i=2,3$, which are of the form
$\overline{\mathcal{S}_{i-1}} \, \E[S_i]/\,
\overline{\mathcal{S}_i}\,\E[S_{i-1}]$ involve ratios of the form
$\E[S_i]/\E[S_{i-1}]$,  which eliminate the $\eta_0$ dependence. A
significant advantage of using ratios
$\overline{\mathcal{S}_i}/\E[S_i]$ is that
$\overline{\mathcal{S}_i}/\E[S_i]=\overline{\varphi_{i}(h_i)}/\phi_i(a_i)$
for $i=1,2$.
 That is, the terms $a_{i}^{2i}$ in the denominators of both
$\overline{\mathcal{S}_i}$ and $\E[S_i]$ drop out -- see equations
(\ref{semi-der2}), (\ref{barS1irr}) for the generalized gradient
constraint, and (\ref{semi-der4}), (\ref{barS2irr}) for the
generalized curvature constraint.
 For example, in the case of the generalized gradient constraint
this means that even if the limit of $\fbarg/a_{1}^{2}$ for $a_1
\rightarrow 0$ is not well defined (i.e., for non-differentiable
models), the ratio $\overline{\mathcal{S}_i}/\E[S_i] \propto \fbarg
/ F_{\la}(a_1)$ is still well defined. Similarly one can show that
the respective ratio for the generalized curvature constraint is
also well defined.

Larger values of the exponent $\beta$ give smaller values of the
distance functional for the same number of iterations. The results
for the SSRF parameters do not depend on $\beta$. Hence, $\beta$ is
a handle on the convergence rate of the optimization and can be set
to one.

Multiple ``solutions'' of the minimization problem for the model
parameters can not be ruled out. The distance functional has by
definition a single solution in terms of the $z_i$.  However, since
the dependence $z_i(\e,\xi,\km)$ is nonlinear, more than one
solutions for $(\e,\xi,\km)$ may be possible, or the optimization
algorithm may get trapped near local minima. It is acceptable to
have more than one solutions corresponding to different types of
``reasonable'' spatial dependence.

\section{Numerical Simulations}
\label{sec:simul}

Numerical experiments based on simulated samples are conducted to
illustrate the performance of the proposed SSRF inference process.
The experiments investigate the ability of FGC-SSRFs to model
spatial distributions generated based on commonly used theoretical
models.
The comparisons are based on the
covariance function and on cross-validation.

\subsection{Covariance Estimation}

\vspace{6pt} Three covariance models are considered:

\begin{enumerate}
\item Spherical
$$ c_s(\bfr)=\sigma_{\rm x}^2 \, \left
\{1-\frac{3}{2}\frac{\|\bfr\|}{b_s}+\frac{1}{2}\frac{\|\bfr\|^3}{b_s^3}
\right\}\bigone_{ 0\leq \|\bfr\| \leq b_s},$$
\item Exponential
$$ c_e(\bfr)= \sigma_{\rm x}^2 \, \exp \left(- \|\bfr\|/b_e \right),$$
\item Gaussian
$$ c_g(\bfr)= \sigma_{\rm x}^2 \, \exp \left(- \|\bfr\|^{2}/b_g^{2} \right).$$
%
%

\end{enumerate}

A uniform distribution of $n=200$ sampling locations $\bfs_1,\ldots
\bfs_n$ on the two-dimensional domain $[0,5]\times [0,5]$ is
assumed. The simulated data are denoted by $X^{(m)}(\bfs_i),$ $1\leq
i \leq n,$ where $1\leq m \leq M,$ is the {\it sample index}. The
data are generated from a standard Gaussian spatial random field
(zero mean and unit variance) using the Cholesky LU decomposition
method. The spatial dependence is given by the three models above,
with correlation lengths $b_s=1$ and $b_e=0.5,$ $b_g=1$. For each
covariance model $M=100$ independent samples are obtained. Each
realization differs from the others in both the sampling locations
and the values of the field. The triangular kernel with support
$[0,1]$ is used in SSRF parameter estimation for all the samples.

\begin{figure}[htbp]
\rotatebox{0}{ \resizebox{10cm}{8cm}{
\includegraphics{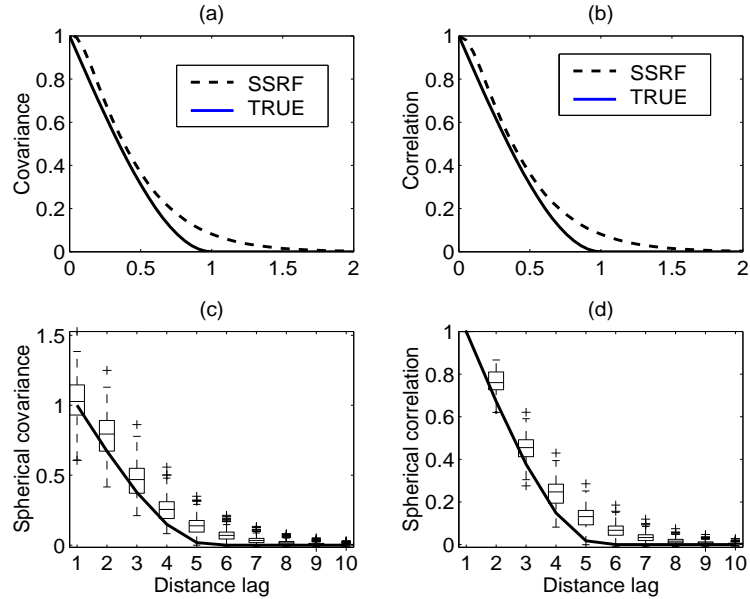}}}
\caption{(a) Spherical covariance (continuous line)  and Spartan
covariance estimator (broken lines). (c) Box plots of Spartan
covariance estimator versus the spherical model. Plots (b) and (d)
are the counterparts  of plots (a) and (c) respectively for the
correlation function.} \label{fig:sphe}
\end{figure}

\begin{figure}[htbp]
\rotatebox{0}{ \resizebox{10cm}{8cm}{
\includegraphics{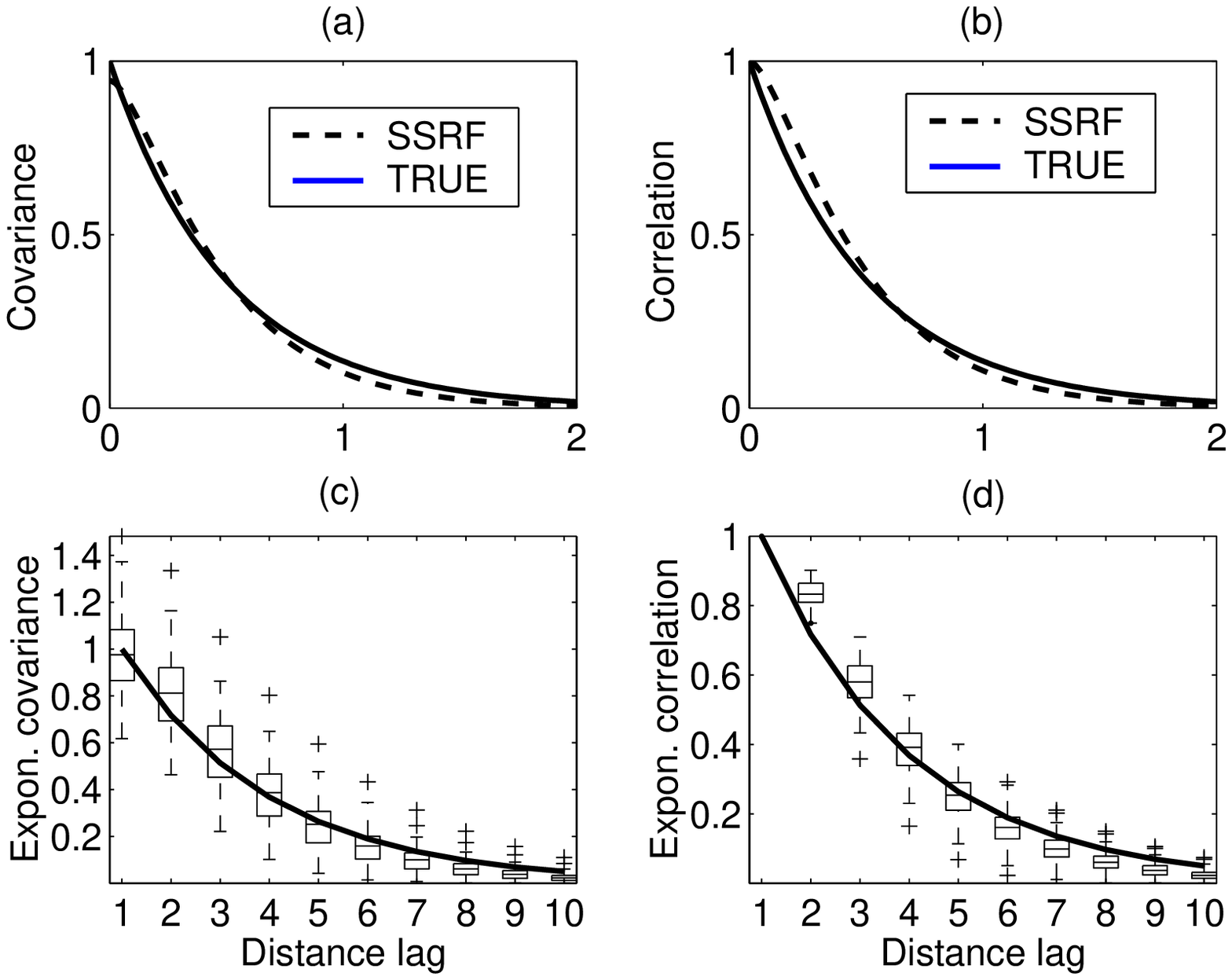}}}
\caption{(a) Exponential covariance (continuous line) and Spartan
covariance estimator (broken lines). (c) Box plots of Spartan
covariance estimator versus the exponential model. Plots (b) and (d)
are the counterparts  of plots (a) and (c) respectively for the
correlation function.} \label{fig:expo}
\end{figure}

\begin{figure}[htbp]
\rotatebox{0}{ \resizebox{10cm}{8cm}{
\includegraphics{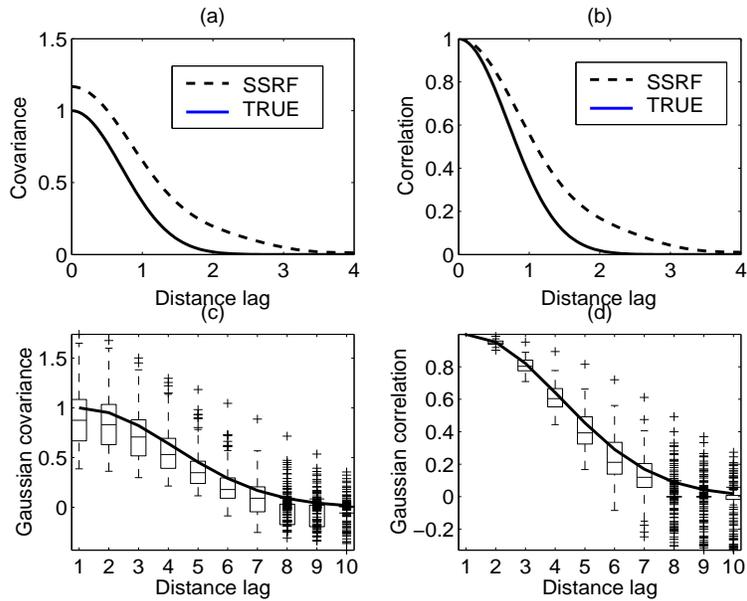}}}
\caption{(a) Gaussian covariance (continuous line)  and Spartan
covariance estimator (broken lines). (c) Box plots of Spartan
covariance estimator versus the Gaussian model. Plots (b) and (d)
are the counterparts of plots (a) and (c) respectively for the
correlation function.} \label{fig:gaus}
\end{figure}

For each sample, the SSRF covariance estimator is calculated at
$10$, uniformly spaced intervals between $0$ and $1.2.$ Figure
\ref{fig:sphe} displays the results for simulations based on the
spherical model. The covariance function obtained from a single
sample is shown for plot (a), and for the correlation function in
plot (b).  The latter is obtained from the SSRF covariance estimator
following division by the sample variance estimate and eliminates
the impact of sample-to-sample variance fluctuations. Box plots
based on all the samples are shown in plot (c) for the covariance
function and in plot (d) for the correlation function. The same
plots for the exponential model are shown in Figure~\ref{fig:expo},
and for the Gaussian model in Figure~\ref{fig:gaus}. The closer
agreement is between the SSRF estimator and the exponential model.
This is justified by the fact that in $d=3$ the SSRF model for $\km
\rightarrow \infty$ and $\eta_1=2$ is equivalent to the exponential
covariance \cite{dth03}.  The SSRF estimator matches the Gaussian
covariance very well near the origin, due to the differentiability
of both models. At large lags the SSRF model box plots exhibit
considerable scatter, which is due to the fact that for certain
realizations the optimization converges to negative $\eta_1$.

It is clear from the plots that the SSRF model does not provide a
perfect match with the ``true'' covariance models over the entire
range of lags. However, this is not a major obstacle in
geostatistical applications, in which the ``true'' covariance is
unknown. In practice, estimation of the empirical covariance
function (or equivalently the variogram) from a single sample
involves considerable uncertainties, which are difficult to quantify
\cite{march}. The uncertainties are more pronounced at larger lags,
at which the averaging procedure involves a smaller number of pairs.
Moreover, the theoretical  covariance functions merely represent
approximations of ``actual'' covariance functions, and thus do not
have any ``inherent'' advantage over the SSRF model.

\subsection{Cross Validation}
\label{ssec:cros-val}

\begin{figure}[htbp]
\rotatebox{0}{ \resizebox{8cm}{6cm}{
\includegraphics{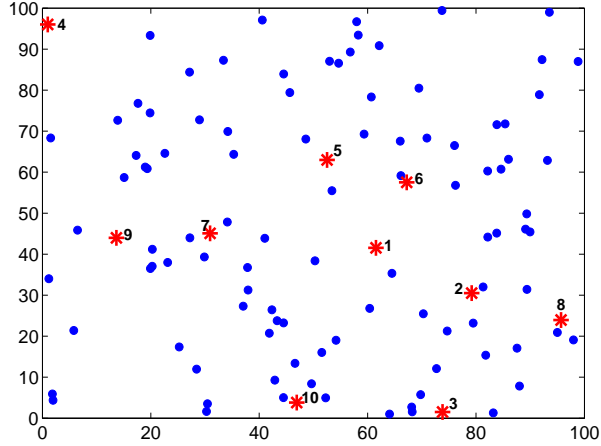}}}
\caption{Locations of the training set (circles) and the validation
set (stars). } \label{fig:locations}
\end{figure}
In geostatistical applications, the performance of a spatial model
is typically evaluated by its ability  to ``predict'' measured
sample values at a number of cross validation points. Here we
consider $n=110$ sampling locations over the domain $D=[0,
100]\times [0,100].$ The set of $110$ points is partitioned into a
validation set, $\mathcal{S}_v$, consisting of $n_{\rm v}=10$ points
chosen at random, and the training set, $\mathcal{S}_t,$ including
the remaining $n_{\rm t}=100$ points. The two sets of points are
shown in Figure~\ref{fig:locations}.

\begin{figure}[htbp]
\rotatebox{0}{ \resizebox{8cm}{6cm}{
\includegraphics{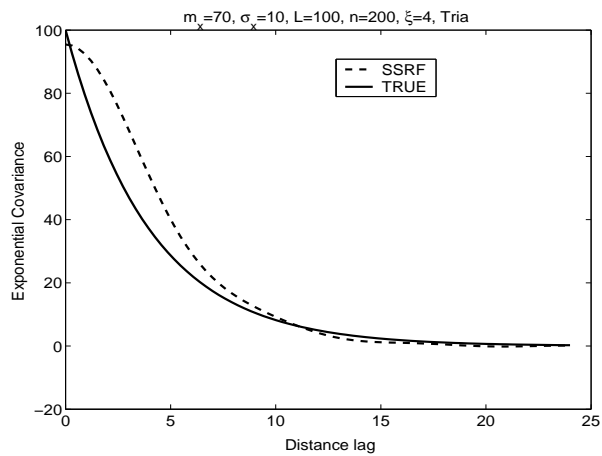}}}
\caption{Comparison of the Spartan covariance estimator (broken
line) with the corresponding theoretical model (continuous line).}
\label{fig:predfig1}
\end{figure}

One hundred independent samples are simulated from a Gaussian SRF
with mean $m_{\rm X}=70$ and standard deviation $\sigma=10$ using an
exponential covariance model, which will henceforth be referred to
as the ``true model''. The correlation length is set to $b_{e}=4$
(i.e., the correlation range, where the covariance drops to $5\%$ of
the initial value is $\approx 12$). The true exponential model and
the SSRF covariance estimator, obtained from a single sample on
$\mathcal{S}_t$, are given in Figure \ref{fig:predfig1}. The
behavior of the Spartan estimator follows the plots of
Figure~(\ref{fig:expo}), that is, the SSRF overestimates the true
model near the origin, where it fails to capture the abrupt decline
of the exponential.

The performance of the SSRF covariance model is evaluated by means
of cross validation.  We use the method of Ordinary Kriging,
e.g.,~\cite{cress}, both with the SSRF covariance and the true
exponential covariance to ``predict'' the field values in
$\mathcal{S}_v.$ The predictions based on the SSRF covariance will
be denoted by $\hat{X}^{(m)}_{\rm ssrf}(\bfs_j),$ while those of the
true model with $\hat{X}^{(m)}_{\rm true}(\bfs_j),$ $\bfs_j \in
\mathcal{S}_v$; $m=1,\ldots,M$ is the realization (sample) index. In
general, a prediction will be denoted by $\hat{X}^{(m)}_{\rm
t}(\bfs_j),$ where ``t=ssrf'' for the SSRF model and ``t=true'' for
the true covariance.  The relative prediction error is then given by
\begin{equation}
\label{eq:pred-re} \epsilon_{{\rm t}}^{(m)}(\bfs_j) =
\frac{\hat{X}^{(m)}_{\rm
t}(\bfs_j)-X^{(m)}(\bfs_j)}{X^{(m)}(\bfs_j)}. \end{equation}

In Table~2 we compare for each point of $\mathcal{S}_v$  the mean
relative error (MRE), $\langle \epsilon_{{\rm t},1}(\bfs_j) \rangle=
\frac{1}{M}\sum_{m=1}^{M} \epsilon_{{\rm t}}^{(m)}(\bfs_j)$, and the
mean absolute relative error (MARE), $\langle \epsilon_{{\rm
t},2}(\bfs_j) \rangle= \frac{1}{M}\sum_{m=1}^{M} |\epsilon_{{\rm
t}}^{(m)}(\bfs_j)|$.  The MRE is $5\%$ or lower for both estimators,
as expected given the fact that kriging is an unbiased predictor.
The MARE is slightly higher for the Spartan model. This is explained
based on the difference between the SSRF and the true covariance
function (see Figure~\ref{fig:predfig1} below). Note that at ${\bf
s}_4$ both models give the same results for the MRE and the MARE.
This happens because ${\bf s}_4$ does not have any nearby neighbors,
and thus the prediction at this point is reduced to the mean value.

The analysis in this section shows that the SSRF covariance model
performs satisfactorily, in terms of cross validation compared to
the predictions obtained with the exponential model used to generate
the data.


\begin{table}
\begin{center}
\begin{tabular}{ccccccccccc}
\hline  \\
  {\bf Statistics}  &  $\bfs_1$ & $\bfs_2$ &
  $\bfs_3$ & $\bfs_4$ &  $\bfs_5$ &$\bfs_6$ &$\bfs_7$ &$\bfs_8$
  &$\bfs_9$ &$\bfs_{10}$
  \\ [0.5ex]
\hline \hline  \\
 MRE True  &  $0.0401 $& $0.0038$ & $-0.0040$ & $0.0615$
 & $0.0533$ & $-0.0126$ & $0.0192 $& $ 0.0630$& $0.0390 $&
 $ 0.0075$ \\
 \\ 
 MRE Spartan & $0.0496 $& $0.0043$ & $-0.0147$ & $0.0615$
 & $0.0535$ & $-0.0102$ & $0.0096 $& $ 0.1177$& $0.0500 $&
 $ 0.0075$ \\
 \hline \hline \\[1ex]
MARE True &  $0.1306 $& $0.0961$ & $0.0944$ & $0.1569$
 & $0.1456$ & $0.0796$ & $0.1052 $& $ 0.1064$& $0.1389 $&
 $ 0.0886$ \\
 \\ 
 MARE Spartan & $0.1319 $& $0.1099$ & $0.1344$ & $0.1569$
 & $0.1502$ & $0.0848$ & $0.1220  $& $ 0.1566$& $0.1552 $&
 $ 0.0956$ \\
 \hline \hline \\
\end{tabular}
\end{center}
\label{tab:pred1} \caption{Prediction errors for  the Spartan
estimator and the true covariance model at the locations in
$\mathcal{S}_v$.  First two rows: MRE. Third and fourth rows: MARE.
The ordering of the validation points is as shown in
Figure~(\ref{fig:locations}).}
\end{table}


\section*{Acknowledgements}
This research is supported by the Marie Curie Action: Marie Curie
Fellowship for the Transfer of Knowledge (Project SPATSTAT, Contract
No. MTKD-CT-2004-014135) and co-funded by the European Social Fund
and National Resources - (EPEAEK-II) PYTHAGORAS.



\begin{thebibliography}{99}

%


\bibitem[Armstrong, 1998]{arm98}
Armstrong, M. (1998). \textit{Basic Linear Geostatistics}. Springer,
Berlin.


\bibitem[Bochner, 1959]{bochner}
Bochner, S. (1959). \textit{Lectures on Fourier Integrals}.
Princeton University Press, Princeton, NJ.


%

\bibitem[Christakos, 1992]{christ}
Christakos, G. (1992). \textit{Random Field Models in Earth
Sciences}. Academic Press, San Diego, CA.

\bibitem[Christakos and Hristopulos, 1998]{ch98}
Christakos, G. and Hristopulos, D. T. (1998). \textit{
Spatiotemporal Environmental Health Modelling}. Kluwer, Boston.

\bibitem[Cressie, 1993]{cress}
Cressie, N. (1993). \textit{Statistics for Spatial Data}. Wiley, New
York.


\bibitem[Gelhar, 1993]{gel93}
Gelhar, L. W. (1993). \textit{Stochastic Subsurface Hydrology}.
Prentice Hall, Englewood Cliffs, NJ.

\bibitem[Goovaerts, 1997]{goov}
Goovaerts, P. (1997). \textit{Geostatistics for Natural Resources
Evaluation}. Oxford University Press, NY.


\bibitem[Hall {\it et al.}, 1994]{hall1}
Hall, P., Fischer, N. and Hoffmann, B. (1994)
\newblock{On the nonparametric
estimation of covariances functions.}  {\em Annals of Statistics},
{\bf 22}, 2115-2134.

\bibitem[Hall and Patil, 1994]{hall2}
Hall, P. and Patil, P. (1994)
\newblock{
 Properties of nonparametric estimators
of autocovariance for stationary random fields}. {\em Probab. Theory
Relat. Fields}, {\bf 99}, 399-424.

\bibitem[Hohn, 1999]{hohn}
Hohn, M. E. (1999). \textit{Geostatistics and Petroleum Geology}.
Kluwer, Dordrecht.

\bibitem[Hristopulos, 2002]{dth02}
Hristopulos, D. T. (2002). New anisotropic covariance models and
estimation of anisotropic parameters based on the covariance tensor
identity. \textit{Stoch. Environ. Res. and Risk Assessment}, {\bf
16}(1), 43-62.

\bibitem[Hristopulos, 2003a]{dth03}
Hristopulos, D. T. (2003). Spartan Gibbs random field models for
geostatistical applications. \textit{SIAM J. Sci. Computing} {\bf
24}, 2125-2162.


\bibitem[Hristopulos, 2004]{dth04}
Hristopulos,  D. T. (2004). Anisotropic Spartan Random Field Models
for Geostatistical Analysis. In \textit{Proceedings of 1st
International Conference on Advances in Mineral Resources Management
and Environmental Geotechnology 2004} (eds Z Agioutantis and K
Komnitsas), pp. 127-132. Heliotopos Conferences, Athens.

\bibitem[Hristopulos, 2005a]{dth05a}
Hristopulos,  D. T. (2005). Identification of Spatial Anisotropy by
means of the Covariance Tensor Identity. In \textit{Mapping
radioactivity in the environment: spatial interpolation comparison
2005} (ed. G. Dubois), Office for Official Publications of the
European Communities, Luxembourg, in print.

\bibitem[Hristopulos, 2005b]{dth05b}
Hristopulos, D. T. (2005).   Spartan gaussian random fields for
geostatistical applications: Non-constrained simulations on square
lattices and irregular grids. \textit{J. Comput. Methods Sci.
Engin.}, {\bf 5}(2), 149-164.


\bibitem[Hristopulos and Elogne, 2006b]{dthse06}
Hristopulos, D. T., and Elogne, S. (2006).
 Analytic Properties and
the covariance function of a class of generalized Gibbs random
fields, cs.IT/0605073 (2006).


\bibitem[Kanevsky and Maignan, 2004]{kan04}
Kanevsky, M. and Maignan, M. (2004). \textit{Analysis and Modelling
of Spatial Environmental Data}. M. Dekker, New York.

\bibitem[Kitanidis, 1997]{kit97}
Kitanidis, P. K. (1997). \textit{Introduction to Geostatistics:
Applications to Hydrogeology}. Cambridge University Press,
Cambridge.

\bibitem[Lantuejoul, 2001]{lantu} Lantuejoul, C. (2001).
\textit{Geostatistical Simulation: Models and Algorithms}. Springer,
Berlin.

\bibitem[Marchant and Lark, 2004]{march}
Marchant, B. P. and Lark, R. M. (2004). Estimating variogram
uncertainty. \textit{Math. Geology}, {\bf 36}(8), 867-898.




\bibitem[Nelder and Mead, 1965]{simplex}
Nelder, J. A. and Mead, R. (1965). A Simplex Method for Function
Minimization. \textit{Comput. J.} {\bf 7}, 308-313.

\bibitem[Rubin, 2003]{rub03}
Rubin, Y. (2003). \textit{Applied Stochastic Hydrogeology}. Oxford
University Press, New York.

%

\bibitem[Smith, 2000]{smith00}
Smith, R.L.  (2000). Spatial statistics in environmental science. In
\textit{Nonlinear and Nonstationary Signal Processing.} (ed. W. J.
Fitzgerald, R. L. Smith, A. T. Walden and P. C. Young), pp. 152-183,
Cambridge University Press, Cambridge.

\bibitem[Stein, 1999]{stein99}
Stein, M.L. (1999). Predicting random fields with increasing dense
observations. \textit{Ann. Appl. Probab.}, {\bf 9}(1), 242-273.






\bibitem[Yaglom, 1987]{yagl87}
Yaglom, A. M. (1987). \textit{Correlation Theory of Stationary and
Related Random Functions I: Basic Results}. Springer, NY.


\end{thebibliography}
 \end{document}